\documentclass[12pt]{article}

\usepackage{amssymb}
\usepackage{amsmath}
\usepackage{graphicx}
\usepackage{hyperref, bookmark}

\newcommand{\BYDEF}{\,\shortstack{\textup{\tiny def} \\ = }\,}
\def\o{\over}

\def\G{\Gamma}
\def\g{\gamma}
\def\gm{\gamma(dy,du)}

\def\a{\alpha}

\def\e{\epsilon}

\def\P{{\mathcal P}}

\def\ph{\varphi}

\newtheorem{Theorem}{Theorem}[section]
\newtheorem{Proposition}[Theorem]{Proposition}

\newtheorem{Lemma}[Theorem]{Lemma}
\newtheorem{Corollary}[Theorem]{Corollary}
\newtheorem{Definition}[Theorem]{Definition}

\numberwithin{equation}{section}

\allowdisplaybreaks

\begin{document}

{\bf \Large{On convergence of occupational measures sets of a discrete-time stochastic control system, with applications to averaging of hybrid systems}}
\bigskip
\bigskip
\bigskip

{\bf L. Gamertsfelder}\\
{\it\small{Department of Mathematics and Statistics, Macquarie University, Eastern Road, Macquarie Park, NSW 2113, Australia}}\\

\bigskip
\bigskip
\bigskip

{\bf Abstract.} \small{In the first part of the paper, we consider a discrete-time stochastic control system. We show that, under certain conditions, the set of random occupational measures generated by the state-control trajectories of the system as well as the set of their mathematical expectations converge (as the time horizon tends to infinity) to a convex and compact (non-random) set, which is shown to coincide with the set of stationary probabilities of the system.  In the second part, we apply the results obtained in the first part to deal with a hybrid system that evolves in continuous time and is subjected to abrupt changes of certain parameters. We show that the solutions of such a hybrid system are approximated by the solutions of a differential inclusion, the right-hand side of which is defined by the limit occupational measures set, the existence and convexity of which is established in the first part of the paper.

\bigskip
\bigskip
\bigskip

\section{Introduction}

In the first part of this paper, we consider a discrete-time stochastic control system. We show that, under certain conditions, the set of random occupational measures generated by the state-control trajectories of the system as well as the set of their mathematical expectations converge (as the time horizon tends to infinity) to a convex and compact (non-random) set, which is shown to coincide with the set of stationary probabilities of the system. 

In the second part of the paper, we use results obtained in the first part to deal with a hybrid system that evolves in continuous time and is subjected to abrupt random changes of certain parameters. The changes of the parameters are synchronized with and controlled by a discrete-time stochastic control system. The controls and states of the latter change their values at the time moments $t_l =l\epsilon, l=1,2,\ldots$ where
$\epsilon$ is a small positive parameter. That is, the systems' parameters change their values frequently (the smaller the $\epsilon$, the higher the frequency). We show that the solutions of such a hybrid system are approximated by the solutions of a differential inclusion, the right-hand side of which is defined by the limit occupational measures set, the existence and convexity of which is established in the first part of the paper. We also indicate the way how a near-optimal control of the hybrid system can be constructed on the basis of an optimal solution of the differential inclusion.

Dynamical systems, in which the state variables evolve in different time scales (so that some of the state variables  can be considered to be slow or fast with respect to others) are often modelled with the help of differential equations that contain
a small  perturbation parameter (or parameters)  near some of the derivatives. This parameter   is  called a singular perturbation parameter (SPP), and the systems containing such a parameter are called singularly perturbed.

Problems of   control  of SP systems appear in multiple applications, and ways of dealing with them
have been discussed extensively in the literature (see  surveys \cite{dmitriev2006singular}, \cite{kokotovic1984applications}, \cite{naidu2002singular}, \cite{zhang2014singular} and also more recent works
\cite{fouque2017perturbation}, \cite{gaitsgory2017near}, \cite{goldys2022singular}, \cite{goreac2010discontinuous}, \cite{yang2020singular} as well as  references therein). One of the most common approaches to SP systems is based on the equating of the SPP to zero.
Being very efficient in dealing with many important classes of optimal control problems
(see, e.g., \cite{bensoussan1988perturbation}, \cite{kabanov2013two}, \cite{kushner2012weak} and \cite{dmitriev2006singular}, \cite{dontchev1996tikhonov}, \cite{kokotovic1999singular}, \cite{veliov1997generalization} for results in the
stochastic and deterministic settings, respectively), it may not
lead to a good approximation of the optimal solution (see examples in \cite{artstein2010periodic}, \cite{borkar2005existence}, \cite{gaitsgory1992suboptimization}, \cite{leizarowitz2002order}).
In fact, the validity of this approach
is related to the validity of the hypothesis that the optimal control  is in some sense slow and that  the fast
state variables converge rapidly to their steady states\footnote{Note that the convergence of optimal state-control trajectories to their steady states is closely related to the ``turnpike phenomenon"; see, e.g.,  \cite{carlson2012infinite}, \cite{zaslavski2014stability}, \cite{zaslavski2014turnpike}.}
and remain in a neighbourhood of these  steady states in the process of the slow evolution. While the validity of such a hypothesis has been established
under natural stability conditions by the celebrated Tikhonov’s theorem in the case of uncontrolled
dynamics (see \cite{tikhonov1952systems}), this hypothesis may not be valid in the control setting, the reason
for this being the fact that the use of rapidly  oscillating  controls may lead to significant
(not tending   with the small parameter to zero) improvements of the performance indices.\\
Various averaging-type approaches allowing one to deal with the state-control trajectories that do not converge to steady states 
have been proposed in the literature (see, e.g., \cite{alvarez2002viscosity}, \cite{artstein1999invariant}, \cite{artstein2000value}, \cite{borkar2005existence}, \cite{borkar2007averaging}, \cite{gaitsgory2004representation}, \cite{gaitsgory1999limit}, \cite{gaitsgory2017near} and references therein). In particular, the averaging technique  based on the convergence of the set of occupational measure generated by the fast state-control trajectories to the set of stationary probability measures for  stochastic and deterministic systems evolving in continuous time were proposed in \cite{borkar2005existence} and, respectively, in  \cite{gaitsgory2004representation}, \cite{gaitsgory1999limit}, and it is this approach that we use to deal with the hybrid system. \\
Hybrid systems similar to the one we are tackling in the paper have been studied in \cite{altman1993control}, \cite{altman1995hybrid},
\cite{altman1997asymptotic}, \cite{brunetti2017nonzero}, and \cite{shi1998asymptotic}, where it was mentioned, in particular, that  these types of systems may arise in the modelling of admission control into telecommunication
networks. 
Note  that, in contrast to the aforementioned works, in which only systems with finite state-control spaces were considered, we deal with a much more general case when the
state-control space   may be non-denumerable, which makes the consideration much more technically
challenging.

The paper is organized as follows. 
In Section \ref{sec:prelim}, we refine the definition of the control plans,
and we introduce some new definitions and notations. In Section \ref{sec:borkar_results_weak}, we introduce some controllability conditions, and we show that, under these controllability conditions,
the set of occupational measures generated by the state-control trajectories with  fixed initial data converges to the set of stationary measures.
In Section \ref{sec:borkar_results_strong}, we introduce stronger controllability conditions, and we show that, under these stronger conditions, the set of random occupational measures generated by the state-control trajectories with  fixed initial data converges to the set of stationary measures. In Section \ref{sec:hybrid_prelim}, we introduce the hybrid system and the optimal control problem of interest. In Section \ref{sec:hybrid_trajectory_approximation}, a differential inclusion is introduced, and it is shown that the slow components of the solutions of the hybrid system are approximated by the solutions of this differential inclusion. In Section \ref{sec:hybrid_system_convergence}, it is shown, under the controllability conditions of Section \ref{sec:borkar_results_strong}, that for every solution of the differential inclusion, there exists a control plan such that the slow component of the corresponding solution of the hybrid system is asymptotically close to this solution of the differential inclusion. Moreover, it is established that an asymptotically optimal control plan for the hybrid control problem can be constructed on the basis of the optimal solution of the differential inclusion. Proofs of the results from Sections \ref{sec:borkar_results_weak} - \ref{sec:borkar_results_strong} are given in
Sections \ref{sec:borkar_proofs_1} - \ref{sec:borkar_proofs_2}.

\section{The Discrete-Time Stochastic Control System and Preliminaries}
\label{sec:prelim}
Consider the   discrete-time stochastic control system in which the state $y(t)$ evolves according to the dynamics
\begin{equation}
	\label{eq:associated_system}
	\begin{aligned}
		y(t+1) = f(y(t), u(t), \sigma(t)), \quad y(0) = y_0, \quad t  = 0, 1, 2, \ldots \ .
	\end{aligned}
\end{equation}
where it is everywhere assumed that the following conditions are satisfied
\begin{itemize}
	\item The function $f(y,u,\sigma) : Y\times \hat{U}  \times S \rightarrow \mathbb{R}^m$ is continuous on $Y\times \hat{U} \times S$, where $Y$ is a compact subset of $\mathbb{R}^m$ and $\hat{U}, S$ are compact metric spaces.
	\item $\sigma(t) \in S$, $t = 0,1,\ldots,$ is a sequence of independent, identically distributed random elements defined on a common probability space.
	\item $f(y,u,\sigma) \in Y$ for any $y \in Y$, any $u \in \hat{U}$,  and any $\sigma \in S$ (that is, the set $Y$ is forward invariant with respect to system (\ref{eq:associated_system})).
	
\end{itemize}

We will be dealing with the special  case when the control set is independent of the state variable. Namely,  $U(y)=\hat U$ for any
$y\in Y$. Note that in this case the graph $G$ of $U(\cdot)$ is given by 
\begin{equation*}\label{e-Y-times-U}
	G=Y\times \hat{U}.
\end{equation*}
Also, the controls that we will be considering are defined as follows.
The control at the moment $t=0$ is chosen as an element of $\hat{U}$, and the control at a moment $t$ ($t=1,2,\ldots$) is chosen as a function of the observed values of $\sigma(0), \ldots, \sigma(t-1)$. That is,
\begin{equation}
	\begin{split}
		\label{eq:control_dependency}
		&u(0) = \pi_0 \in \hat{U}\\
		&u(t) = \pi_t(\sigma(0),  \ldots, \sigma(t-1)) \in \hat{U}.
	\end{split}
\end{equation}
where the functions $\pi_t(\cdot),$  $t=1,2,\ldots$ are assumed to be Borel measurable. We will refer to any sequence of such functions $\pi = \{\pi_t(\cdot), \ t = 0, 1, \ldots \}$ as a \textit{control plan}, and we will be denoting
the set of all controls plans as $\Pi$.






Given a control plan $\pi$ and an initial condition $y(0)=y_0\in Y$,\ let $(y^{\pi, y_0}(\cdot),u^{\pi, y_0}(\cdot))$  stand  for the state-control trajectory of system (\ref{eq:associated_system}) and let $\g^{\pi , y_0 ,T}\in \mathcal{P}(Y \times \hat{U}) $ be the occupational measure generated by $\pi$ over the time sequence
$ \{ 0,1,...,T-1 \}$:
\begin{equation}\label{E6-exp-taime-ave}
	\g^{\pi, y_0, T}(Q)={1\over T} E\left[\sum_{t=0}^{T-1} 1_Q(y^{\pi, y_0}(t),u^{\pi, y_0}(t))\right]\ \ \ {\rm for \ any \ Borel}\ \ \  Q\subset Y \times \hat{U} .
\end{equation}
Denote  by $ \G_T(y_0)$ the set of all such occupational measures:
\begin{equation}\label{union-2-occ-1}
	\G_T(y_0)\BYDEF\bigcup_{\pi \in \Pi}\{\g^{\pi, y_0, T}\},
\end{equation}
the union being over all control plans.

To describe convergence properties of occupational measures,  the following metric on
$\P(G)$ will be used:
\begin{equation}\label{e-rho}
	\rho(\g',\g''):=\sum_{j=1}^{\infty} {1\o 2^j}\left|\int_G q_j(y,u)\g'(dy,du)-\int_G q_j(y,u)\g''(dy,du)\right|
\end{equation}
for $\g',\g''\in \P(G)$, where $q_j(\cdot),\,j=1,2,\dots,$ is a sequence of Lipschitz continuous functions dense in the unit ball of the space of continuous functions $C(G)$ from $G$ to $\mathbb{R}$.
This metric is consistent with the weak$^*$ convergence topology on $\P(G)$, that is,
a sequence $\g^k\in \P(G)$ converges to $\g\in \P(G)$ in this metric if and only if
$$
\lim_{k\to \infty}\int_G q(y,u)\g^k(dy,du)=\int_G q(y,u)\g(dy,du)
$$
for any $q\in C(G)$.

The Hausdorff metric $\rho_H(\Gamma_1,\Gamma_2)$ between $\Gamma_1\subset \P(G)$ and $\Gamma_2\subset \P(G)$ is defined as follows:

\begin{equation}\label{eq:measures-distance-definitions}
	\rho(\g,\Gamma)\BYDEF\inf_{\g'\in \Gamma}\rho(\g,\g'),\quad
	\rho_H(\Gamma_1,\Gamma_2)\BYDEF\max\{\sup_{\g\in \Gamma_1}\rho(\g,\Gamma_2),\sup_{\g\in \Gamma_2}\rho(\g,\Gamma_1)\}.
\end{equation}

\begin{Definition}
	A convex and compact set $\G \subset \mathcal{P}(Y \times \hat{U})$ will be called the \textit{limit occupational measures set} (LOMS) of the system (\ref{eq:associated_system}) if
	\begin{equation}\label{eq:LOMS}
		\rho_H(\G_T(y_0), \G) \le \nu^*(T) \ \ \ \forall \ y_0\in Y, \ \ \ \ {\rm where }\ \ \ \lim_{T \rightarrow \infty} \nu^*(T) = 0.
	\end{equation}
\end{Definition}

Define the set $ W\subset \P(G)$ by the equation
\begin{equation}\label{eq:s_stationary}
	\begin{aligned}
		W\BYDEF\left\{\g\in \P(G) |\, \int_{G}\big(E\big[\ph(f(y,u,\sigma))\big]-\ph(y)\big)
		\gm=0\ \ \forall \ph\in C(Y)\right\}.
	\end{aligned}
\end{equation}
Note that, as can be readily seen, the set $W$ is convex and weak$^*$ compact in the metric $\rho$ defined in (\ref{eq:measures-distance-definitions}). For convenience, we refer to the set $W$ as the set of stationary measures.

\begin{Lemma}\label{prop:G_equals_W}
	If the LOMS $\Gamma$ of the system (\ref{eq:associated_system}) exists, then
	\begin{equation}\label{eq:G_equals_W}
		\G = W,
	\end{equation}
	where $W$ is the set of stationary measures defined in (\ref{eq:s_stationary}).
\end{Lemma}
{\bf Proof.}
Let $\G$ be the LOMS of the system (\ref{eq:associated_system}). From (\ref{eq:LOMS}) and the fact that $\G$ is convex, it follows that
\begin{equation}\label{e-rho-H-extra}
	\rho_H(\bar{\rm co}\G_{T},\Gamma) \leq \rho_H(\G_{T},\Gamma)
	\leq \nu^*(T),
\end{equation}
where $\G_T\BYDEF\bigcup_{y_0\in Y}\{\G_T(y_0)\}$.
From (\ref{e-rho-H-extra}) and Theorem 3.2 in \cite{avrachenkov2022lp},  it follows that
\begin{equation*}
	\rho_H(\G, W) \le \lim_{T \rightarrow \infty}\rho_H(\G, \bar{\rm co}\G_T) + \lim_{T \rightarrow \infty}\rho_H(\bar{\rm co}\G_T, W) = 0,
\end{equation*}
which proves (\ref{eq:G_equals_W}).
$\ \Box$

Along with the occupational measure   $\g^{\pi , y_0 ,T}$  defined by (\ref{E6-exp-taime-ave}), let us define the {\it random occupational measure}
$\chi^{\pi , y_0 ,T}\in \mathcal{P}(Y \times \hat{U}) $ generated by a control plan $\pi$ (with an initial condition $y(0)=y_0$) by the equation
\begin{equation}\label{eq:random-occ-measure}
	\chi^{\pi, y_0, T}(Q)={1\over T} \sum_{t=0}^{T-1} 1_Q(y^{\pi, y_0}(t),u^{\pi, y_0}(t)) \ \ \ {\rm for \ any \ Borel}\ \ \  Q\subset Y \times \hat{U} .
\end{equation}
From (\ref{eq:random-occ-measure}) it follows that
\begin{equation}\label{G88}
	\int_{Y \times \hat{U}} q(y,u) \chi^{\pi, y_0, T}(dy,du)={1\over T}\sum_{t=0}^{T-1} q(y^{\pi, y_0}(t),u^{\pi, y_0}(t)),
\end{equation}
for any bounded Borel measurable function $q$ on $G$.


Denote by  $B_T(y_0) $  the set of all random occupational measures
\begin{equation}\label{union-2-occ-2}
	B_T(y_0)\BYDEF\bigcup_{\pi \in \Pi}\{\chi^{\pi, y_0, T}\},
\end{equation}
where the union is over all control plans. We will be interested in investigating  the convergence of the set $B_T(y_0) $ to the set of stationary measures $W$ with $T$ tending to infinity. To characterise such a convergence, we will use the following definition of the \lq\lq  distance" $\rho^E(\g,\Phi)$  between a given (\lq\lq deterministic") probability measure $\g \in \P(Y \times \hat{U})$ and a collection  $\Phi\subset \P(Y \times \hat{U}) $ of random probability measures:
\begin{equation}\label{eq:expected_rho_set_distances}
	\rho^E(\g,\Phi) \BYDEF \inf\limits_{\g' \in \Phi}E[\rho(\g,\g')],
\end{equation}
where  $\rho(\cdot,\cdot)$ is as defined  in (\ref{e-rho}). In addition, we will be using the following Hausdorff type  distance between collections of random elements taking values in $\mathbb{R}^j$. Let $V_1$ and $V_2$ be two collections of integrable random variables defined on the same probability space and taking values in $\mathbb{R}^j$. We define the distance between a random vector $v \in \mathbb{R}^j$ and $V_i$ by the equation
\begin{equation}\label{def:expected_distance}
	d^E(v,V_i) \BYDEF \inf\limits_{v' \in V_i} E[||v - v'||],\ \ \ i=1,2,
\end{equation}
and we define the Hausdorff type distance between $V_1$ and $V_2$ by the equation
\begin{equation}\label{def:expected_hausdorff_distance}
	d^E_H(V_1, V_2) \BYDEF \max\left\{\sup_{\zeta \in V_1}d^E(\zeta, V_2), \sup_{\zeta \in V_2}d^E(\zeta, V_1) \right\}.
\end{equation}
Note that $d^E_H(V_1, V_2)=d_H(V_1, V_2) $ if $V_1$ and $V_2$ are \lq\lq non-random" bounded subsets of $\mathbb{R}^j$, where
$d_H(V_1, V_2)$ is the \lq\lq standard" Hausdorff metric:
\begin{equation}\label{def:hausdorff_distance}
	d_H(V_1, V_2) = \max\left\{\sup_{\zeta \in V_1}d(\zeta, V_2), \sup_{\zeta \in V_2}d(\zeta, V_1) \right\},
\end{equation}
with
\begin{equation}\label{def:distance}
	d(v,V_i) = \inf\limits_{v' \in V_i} ||v - v'||, \ \ \ i=1,2.
\end{equation}

\section{Convergence of the set of occupational measures to the set of stationary measures}
\label{sec:borkar_results_weak}

Let $h(y,u) : Y \times \hat{U} \rightarrow \mathbb{R}^j$ be defined as follows 
\begin{equation}
	\label{eq:h-AC_function}
	h(y,u) \BYDEF (q_1(y,u),q_2(y,u),\ldots,q_j(y,u)),
\end{equation}
where $q_i(y,u),$ $i=1,2,\ldots,$ is a sequence of Lipschitz continuous functions that is dense in the unit ball of the space $C(Y \times \hat{U})$.
\begin{Definition}[weak $h$-approximation conditions]
	\label{def:W-h-AC}
	The weak $h$-approximation conditions are said to be satisfied for system (\ref{eq:associated_system}) if, corresponding to any initial condition $y_0' \in Y$, any control plan $\pi'$ and any other initial condition $y_0''\in Y$, there exists a control plan $\pi''$ such that the state-control trajectories obtained with the use of $\pi'$ and $\pi''$ in system (\ref{eq:associated_system}) satisfy the inequality
	\begin{equation}
		\left|\left|E\left[\frac{1}{T}\sum_{t=0}^{T}h(y^{\pi',y_0'}(t), u^{\pi',y_0'}(t))\right] - E\left[\frac{1}{T}\sum_{t=0}^{T}h(y^{\pi'', y_0''}(t), u^{\pi'', y_0''}(t))\right]\right|\right| \le \nu_h(T),
	\end{equation}
	for some monotone decreasing $\nu_h(\cdot,\cdot):[0,\infty) \rightarrow [0,\infty)$ such that $\lim\limits_{T\rightarrow \infty}\nu_h(T) = 0$.
\end{Definition}
Denote by $V_h(T,y_0)$ the following subset of $\mathbb{R}^j$:
\begin{equation}\label{eq:V_h_T_definition}
	V_h(T,y_0) \BYDEF \bigcup_{\g \in \G_T(y_0)}\left(\int_{Y \times \hat{U}}^{} h(y ,u) \g(dy, du) \right)
\end{equation}
where $h(\cdot, \cdot)$ is as in (\ref{eq:h-AC_function}).

\begin{Theorem}\label{thm:expected_occ_set_and_V_h_distance}
	If the system (\ref{eq:associated_system}) satisfies the weak $h$-approximation conditions (see Definition \ref{def:W-h-AC}) for any vector function $h(\cdot,\cdot)$ as in (\ref{eq:h-AC_function}), then there exists a convex and compact set $V_h \in \mathbb{R}^j$ such that, for any initial condition $y_0\in Y$, \begin{equation}\label{eq:expected_occ_set_and_V_h_distance}
		d_H(V_h(T,y_0), V_h) \le \nu'_h(T), \ \ \ \ \lim_{T \rightarrow \infty}\nu'_h(T) = 0.
	\end{equation}
	Conversely, if there exists $V_h$ such that (\ref{eq:expected_occ_set_and_V_h_distance}) is valid for any initial condition $y_0$, then the weak $h$-approximation conditions are satisfied for any initial conditions $y'_0$, $y''_0$.
\end{Theorem}
{\bf Proof.}
The proof is given in Section \ref{sec:borkar_proofs_1}. $\ \Box$

\begin{Theorem}\label{thm:LOMS_theorem_2}
	If the system (\ref{eq:associated_system}) satisfies the weak $h$-approximation conditions (see Definition \ref{def:W-h-AC}) for any vector function $h(\cdot,\cdot)$ as in (\ref{eq:h-AC_function}), then the LOMS $\G$ of the system (\ref{eq:associated_system}) exists and allows the following representation:
	\begin{equation}\label{eq:LOMS_theorem_2_2}
		\G \BYDEF \{\g \in \mathcal{P}(Y \times \hat{U}) | \int_{Y \times \hat{U}}^{} h(y,u) \g(dy,du) \in V_h \ \ \forall \ h(y,u) \ \text{as in} \ (\ref{eq:h-AC_function})  \},
	\end{equation}
	where $V_h$ are convex and compact sets the existence of which (for every $h(\cdot,\cdot)$ as in (\ref{eq:h-AC_function})) is established by Theorem \ref{thm:expected_occ_set_and_V_h_distance}. That is, for any initial condition $y_0 \in Y$, the following estimate is valid:
	\begin{equation}\label{eq:LOMS_theorem_2_1}
		\rho_H(\G_T(y_0), \G) \le \nu^{*}(T), \ \ \ \ \lim_{T \rightarrow \infty}\nu^{*}(T) = 0,
	\end{equation}
	where $\G$ is as in (\ref{eq:LOMS_theorem_2_2}).
	Conversely, if there exists a convex and compact set $\G \in \mathcal{P}(Y \times \hat{U})$ which satisfies (\ref{eq:LOMS_theorem_2_1}) for any initial condition $y_0$, then the weak $h$-approximation is satisfied for any vector function $h(\cdot, \cdot)$ as in (\ref{eq:h-AC_function}) and any initial conditions $y'_0$, $y''_0$. Also, for any $h(\cdot, \cdot)$ as in (\ref{eq:h-AC_function}), the estimate (\ref{eq:expected_occ_set_and_V_h_distance}) is valid with
	\begin{equation}\label{eq:V_h_G_representation}
		V_h = \bigcup_{\g \in \G} \left\{\int_{Y \times \hat{U}}^{} h(y,u) \g (dy,du) \right\}.
	\end{equation}
	
\end{Theorem}
{\bf Proof.}
The proof is given in Section \ref{sec:borkar_proofs_1}. $\ \Box$

\begin{Corollary}\label{cor:LOMS_theorem_2}
	If the system (\ref{eq:associated_system}) satisfies the weak $h$-approximation conditions (see Definition \ref{def:W-h-AC}) for any vector function $h(\cdot,\cdot)$ as in (\ref{eq:h-AC_function}), then \textnormal{(i)} the set $\G$ defined in (\ref{eq:LOMS_theorem_2_2}) satisfies the equality $\G = W$ and \textnormal{(ii)} the following inequality is valid
	\begin{equation}\label{eq:LOMS_theorem_2_1_W}
		\rho_H(\G_T(y_0), W) \le \nu^{*}(T), \ \ \ \ \lim_{T \rightarrow \infty}\nu^{*}(T) = 0.
	\end{equation}
\end{Corollary}
{\bf Proof.}
Part \textnormal{(i)} follows from Lemma \ref{prop:G_equals_W} and the fact that $\G$ is the LOMS of the system (\ref{eq:associated_system}). Part \textnormal{(ii)} follows from \textnormal{(i)} and (\ref{eq:LOMS_theorem_2_1}). $\ \Box$

Note that a similar inequality to (\ref{eq:LOMS_theorem_2_1_W}) has been proven in Theorem 3.2 of \cite{avrachenkov2022lp} without any controllability assumptions. However, the convergence of the result given there is for the convex hull of the union of $\G_T(y_0)$ over all  initial conditions in $Y$, while the convergence of $\G_T(y_0)$ to $W$ in (\ref{eq:LOMS_theorem_2_1_W})  is for an \lq\lq individual" initial condition. That is, the statement about the validity of  (\ref{eq:LOMS_theorem_2_1_W}) is much stronger. 

\section{Convergence of the set of random occupational measures to the set of stationary measures}
\label{sec:borkar_results_strong}
Let $h(y,u) : Y \times \hat{U} \rightarrow \mathbb{R}^j$ be the function defined in Section \ref{sec:borkar_results_weak}. 
\begin{Definition}[strong $h$-approximation conditions]
	\label{def:S-h-AC}
	The strong $h$-approximation conditions are said to be satisfied for system (\ref{eq:associated_system}) if, corresponding to any initial condition $y_0' \in Y$, any control plan $\pi'$ and any other initial condition $y_0''\in Y$, there exists a control plan $\pi''$ such that the state-control trajectories obtained with the use of $\pi'$ and $\pi''$ in system (\ref{eq:associated_system}) satisfy the inequality
	\begin{equation}
		E\left[\left|\left|\frac{1}{T}\sum_{t=0}^{T}h(y^{\pi',y_0'}(t), u^{\pi',y_0'}(t)) - \frac{1}{T}\sum_{t=0}^{T}h(y^{\pi'', y_0''}(t), u^{\pi'', y_0''}(t))\right|\right|\right] \le \nu_h(T)
	\end{equation}
	for some monotone decreasing $\nu_h(\cdot,\cdot):[0,\infty) \rightarrow [0,\infty)$ such that $\lim\limits_{T\rightarrow \infty}\nu_h(T) = 0$.
\end{Definition}
Note that the weak $h$-approximation conditions are satisfied for the system (\ref{eq:associated_system}) whenever the strong $h$-approximation conditions are satisfied for the system (\ref{eq:associated_system}). In particular, if the strong $h$-approximation conditions are satisfied for the system (\ref{eq:associated_system}), then the results from Section \ref{sec:borkar_results_weak} also hold. Note also that, as demonstrated by the example below, the strong $h$-approximation conditions are satisfied when the control system is linear.

{\bf Example 3.} Consider the system described by the equation
\begin{equation}\label{eq:strong-h-example-system}
	y(t+1) = Ay(t) + B(\sigma(t))u(t) + \sigma(t), \ \ t = 0,1,2,\ldots,
\end{equation}
where $A$ is an $m\times m$ matrix such that $\|A \| \le 1$. 
Let $\pi$ be an arbitrary control plan and let $y_0'$,  $y_0'' \in Y$. We have
$$
y^{\pi,y_0'}(t)-y^{\pi,y_0''}(t) = A (y^{\pi,y_0'}(t-1)-y^{\pi,y_0''}(t-1)) 
$$
$$
\Rightarrow \| y^{\pi, y_0'}(t) - y^{\pi, y_0''}(t)\| \le \| A \|^t \| y_0' - y_0'' \|. 
$$
Then
\begin{equation*}
	\begin{split}
		E\left[\left|\left|\frac{1}{T}\sum_{t=0}^{T}h(y^{\pi,y_0'}(t), u^{\pi,y_0'}(t)) - \frac{1}{T}\sum_{t=0}^{T}h(y^{\pi,y_0''}(t), u^{\pi,y_0''}(t))\right|\right|\right] &\leq \frac{1}{T}\sum_{t=0}^T L_h ||y_0'-y_0''|| ||A||^t \\
		& \leq \frac{1}{T}L_h ||y_0'-y_0''||\frac{1}{1-||A||},
	\end{split}
\end{equation*}
where $L_h \ge 0$ is the Lipschitz constant of the function $h(\cdot,u)$. Thus, the strong $h$-approximation conditions are  satisfied for the system (\ref{eq:strong-h-example-system}), with $\pi'' = \pi'$ and $\nu_h(T) =\frac{1}{T}L_h\frac{1}{1-||A||}\max\limits_{y_0',y_0''\in Y} ||y_0'-y_0''||$.

\begin{Theorem}
	\label{thm:key_theorem_from_borkar}
	If the system (\ref{eq:associated_system}) satisfies the strong $h$-approximation conditions (see Definition \ref{def:S-h-AC}) for any vector function $h(\cdot,\cdot)$ as in (\ref{eq:h-AC_function}), then, 
	\begin{equation}\label{eq:main_borkar_result}
		\sup_{\g \in \G} \rho^E(\g,B_T(y_0)) \le \bar{\nu}(T), \ \ \ \ \lim_{T \rightarrow \infty}\bar{\nu}(T) = 0,
	\end{equation}
	for any initial condition $y_0 \in Y$, where $\G$ is as defined in (\ref{eq:LOMS_theorem_2_2}) and where $\rho^E(\cdot, \cdot)$ is as defined in (\ref{eq:expected_rho_set_distances}).
\end{Theorem}
{\bf Proof.}
The proof is given in Section \ref{sec:borkar_proofs_1}. $\ \Box$

\begin{Corollary}\label{cor:key_theorem_from_borkar}
	If the system (\ref{eq:associated_system}) satisfies the strong $h$-approximation conditions (see Definition \ref{def:S-h-AC}) for any vector function $h(\cdot,\cdot)$ as in (\ref{eq:h-AC_function}), then, 
	\begin{equation}\label{eq:main_borkar_result_W}
		\sup_{\g \in W} \rho^E(\g,B_T(y_0)) \le \bar{\nu}(T), \ \ \ \ \lim_{T \rightarrow \infty}\bar{\nu}(T) = 0,
	\end{equation}
	for any initial condition $y_0 \in Y$.
\end{Corollary}
{\bf Proof.} Follows readily from Theorem \ref{thm:key_theorem_from_borkar} and Corollary \ref{cor:LOMS_theorem_2}\textnormal{(i)}. $\ \Box$


\begin{Theorem}\label{thm:G_T_W_distance}
	There exists a function $\nu^*(\cdot)$, $\lim\limits_{T \rightarrow \infty} \nu^*(T) = 0$, such that 
	\begin{equation}\label{eq:G_T_W_distance}
		\sup_{\chi \in B_T(y_0)} E[\rho(\chi,W)] \le \nu^*(T)
	\end{equation}
	for any initial condition $y_0 \in Y$, where $\rho(\chi, W)$ is as defined in (\ref{eq:measures-distance-definitions}).
\end{Theorem}
{\bf Proof.} The proof is given in Section \ref{sec:borkar_proofs_2}. $\ \Box$

From (\ref{eq:main_borkar_result}) and (\ref{eq:G_T_W_distance}), we have that, under the strong $h$-approximation conditions,  the set of random occupational measures generated by the state-control trajectories
converges to the set of stationary measures $W$.

To conclude this section, we provide the following pair of corollaries which will be utilised when averaging the hybrid system introduced in the second part of the paper.

Denote by $V_h^r(T,y_0)$ the following subset of $\mathbb{R}^j$:
$$V_h^r(T,y_0) \BYDEF \bigcup_{\chi \in B_T(y_0)}\left(\int_{Y \times \hat{U}}^{} h(y ,u) \chi(dy, du) \right).$$
\begin{Corollary}
	\label{cor:vg_vgsy0_distance}
	If the system (\ref{eq:associated_system}) satisfies the strong $h$-approximation conditions for any vector function $h(\cdot,\cdot)$ as in (\ref{eq:h-AC_function}), then the limit set $V_h$ (see (\ref{eq:expected_occ_set_and_V_h_distance})) is presentable in the form: 
	\begin{equation}\label{eq:V_h_W_representation}
		V_h = \bigcup_{\g \in W} \left\{\int_{Y \times \hat{U}}^{} h(y,u) \g (dy,du) \right\}.
	\end{equation}
	Also
	\begin{equation}\label{eq:main_borkar_result_d_E}
		\sup\limits_{v\in V_h}d^E(v,V_h^r(T,y_0)) \le \bar{\nu}_h(T), \ \ \ \ \lim_{T \rightarrow \infty}\bar{\nu}_h(T) = 0,
	\end{equation}
	for any initial condition $y_0 \in Y$, where $d^E(v, V_h^r(T,y_0))$ is as defined in (\ref{def:expected_distance}), and $V_h$ is as in (\ref{eq:V_h_G_representation}).
\end{Corollary}
{\bf Proof.}
The equality (\ref{eq:V_h_W_representation}) follows from (\ref{eq:V_h_G_representation}) and Corollary \ref{cor:LOMS_theorem_2}\textnormal{(i)}. The inequality  (\ref{eq:main_borkar_result_d_E}) follows from (\ref{eq:main_borkar_result_W}) and (\ref{eq:V_h_W_representation}).$\ \Box$

\begin{Corollary}
	\label{cor:vg_vgsy0_distance_expectation}
	There exists a function $\nu^*_h(\cdot)$,  $\lim\limits_{T\rightarrow \infty}\nu^*_h(T) =0$, such that 
	\begin{equation}\label{eq:Vh_VhTy0_distance}
		\sup_{v \in V_h^r(T,y_0)}E[d(v,V_h)] \le \nu^*_h(T) 
	\end{equation}
	for any initial condition $y_0 \in Y$.
\end{Corollary}
{\bf Proof.} 
The proof follows directly from Theorem \ref{thm:G_T_W_distance}. $\ \Box$

\section{The Hybrid Control System and Preliminaries}
\label{sec:hybrid_prelim}
Consider the stochastic control system in which the continuously changing vector of the state variables $z_\epsilon(s) \in \mathbb{R}^n$ evolves according to the equations:
\begin{subequations}
	\label{eq:hybrid_system}
	\begin{align}\label{eq:hybrid_system_1}
		&\dot{z}_\epsilon(s) = g(z_\e(s), y_\e(s), u_\epsilon(s)) \ \ \ \forall \ s \in [0, 1],   \\\label{eq:hybrid_system_2}
		&y_\e(s) = y(t),\ u_\epsilon(s) = u(t) \ \ \ \ \ \forall \ s \in [\epsilon t, \epsilon(t+1)), \ t=0,1,\ldots,\lfloor 1 / \e \rfloor,\\\label{eq:hybrid_system_3}
		&z_\e(0) = z_{ 0},
	\end{align}
\end{subequations}
where $y(t)$ and $u(t)$ are defined in accordance with (\ref{eq:associated_system}) - (\ref{eq:control_dependency}), $\e > 0$ is a small parameter, and $\lfloor \cdot \rfloor$ stands for the floor function. In this chapter we will be dealing with the control plans $\pi$ defined similarly to (\ref{eq:control_dependency}) on the sequence of time moments $t=0,1, \ldots, \lfloor 1 / \epsilon \rfloor$. The set of all such control plans will be denoted as $\Pi$. Given a control plan $\pi\in \Pi$ and the pair $(y^{\pi, y_0}(\cdot), u^{\pi, y_0}(\cdot)) $ obtained in accordance with (\ref{eq:associated_system}) - (\ref{eq:control_dependency}),
we denote by $(y_{\epsilon}^{\pi, y_0}(\cdot), u_{\epsilon}^{\pi, y_0}(\cdot))$ the pair of functions obtained by substituting $(y^{\pi, y_0}(\cdot), u^{\pi, y_0}(\cdot))$ into (\ref{eq:hybrid_system_2}), and we denote by $z^{\pi, y_0}_{\epsilon}(\cdot) $ the respective solution obtained from (\ref{eq:hybrid_system_1}).

The controls and the states of the underlying discrete-time parameters of the system (\ref{eq:hybrid_system}) change their values at the time moments $ \epsilon t, \   t = 0, 1,\ldots \lfloor 1 / \epsilon \rfloor .$   That is, the systems' parameters change their values frequently (the smaller is
$\epsilon$, the higher is the frequency), and,  while the change in values of these parameters from a moment $\epsilon t $ to the moment $\epsilon (t+1) $ is finite (not tending with $\epsilon$ to zero),
the change of the continually evolving state variables on the interval $[\epsilon t  , \epsilon(t+1)] $ is of the order $O(\epsilon) $, this allowing one to characterise the former as \lq\lq fast" and the latter as \lq\lq slow" components of the dynamics.

The function $g(\cdot):\mathbb{R}^n\times Y \times \hat{U}\rightarrow \mathbb{R}^n$ in (\ref{eq:hybrid_system_1}) is assumed to be continuous and satisfy Lipschitz conditions in $z$ on a sufficiently large compact subset $Z$ of $\mathbb{R}^n$ that contains all the solutions of (\ref{eq:hybrid_system_1}) in its interior. That is, for some positive $C > 0$,
\begin{equation}
	\label{eq:Lipschitz_g_assumption}
	\left|\left|g(z_1, y, u) - g(z_2, y, u)\right|\right| \le C||z_1 - z_2|| \quad \forall \  z_1, z_2 \in Z, \ \forall \  y \in Y, \ \forall \ u \in \hat U,
\end{equation}
Note that, due to continuity of $g(\cdot)$, 	
\begin{equation}
	\label{eq:Lipschitz_bounded_assumption}
	\left|\left|g(z,y,u)\right|\right| \le M \ \ \forall \ z \in Z, \ \forall \ y \in Y,  \ \forall \ u \in \hat U,
\end{equation}
where $M$ is a positive constant. 
%


On the solutions of system (\ref{eq:hybrid_system}), we consider the following optimal control problem
\begin{equation}
	\label{eq:OC_hybrid}
	F^\epsilon\BYDEF \inf_{\pi \in \Pi} E[G(z^{\pi, y_0}_\e(1))]
\end{equation}
where $G(\cdot): Z \rightarrow \mathbb{R}$ is Lipschitz continuous with Lipschitz constant $C_G$ and the infimum is sought over the set of control plans from $\Pi$. 

\section{Approximation of the $z$-components of the hybrid system trajectories by the solutions of a differential inclusion}
\label{sec:hybrid_trajectory_approximation}

Let $V_g(z)$ be defined by the equation:
\begin{equation}
	\label{def:vg}
	V_g(z) \BYDEF \bigcup_{\g \in W}\left\{\int_{Y \times \hat{U}}^{}g(z,y,u)\g(dy,du)\right\}.
\end{equation}
where $g(\cdot)$ is as in (\ref{eq:hybrid_system}). Note that 
from the fact that the set $W$ is convex and weak$^*$ compact, it follows that the set $V_g(z)$ is convex and compact. Also, from (\ref{eq:Lipschitz_g_assumption}) and (\ref{eq:Lipschitz_bounded_assumption}), it follows that
\begin{equation}
	\label{eq:V_g_lipschitz_continuity}
	d_H(V_g(z'),V_g(z'')) \le C||z' - z''|| \quad \forall \ z', z'' \in Z
\end{equation}
and
\begin{equation}
	\label{eq:V_g_bounded}
	\max\left\{||v|| : v \in V_g(z) \right\} \le M.
\end{equation}
Consider the differential inclusion
\begin{equation}
	\label{eq:diff_inclusion}
	\dot{z}(s) \in V_g(z(s)), \quad z(0) = z_0
\end{equation}
with the associated optimal control problem
\begin{equation}\label{eq:OC_diff_inclusion}
	F^0 \BYDEF \inf_{\{z(\cdot)\}} G(z(1)),
\end{equation}
where the infimum is on the set of the solutions of the differential inclusion (\ref{eq:diff_inclusion}). The approximation of the $z$-components of the system (\ref{eq:hybrid_system}) by solutions of the above differential inclusion is established by the following theorem.

\begin{Theorem}
	\label{thm:first_hybrid_theorem}
	Let $\pi$ be an arbitrary control plan and let $z^{\pi, y_0}_{\e}(\cdot)$ be the corresponding solution of (\ref{eq:hybrid_system}). There exists  $z_\e(\cdot)$, a solution of the differential inclusion (\ref{eq:diff_inclusion}), such that 
	\begin{equation}\label{eq:DI_sol_existence}
		\max_{s \in [0, 1]} E || z^{\pi, y_0}_\e(s) - z_\e(s)|| \le \psi_1(\epsilon)
	\end{equation}
	for some $\psi_1(\epsilon)$ such that $\lim\limits_{\e \rightarrow 0}\psi_1(\epsilon) = 0$.
\end{Theorem}
{\bf Proof.} The proof of this theorem is given at the end of this section. It is based on Lemmas \ref{lem:K_t} - \ref{lem:gronwall} that are stated and
proved below. $\ \Box$\\
Choose a function $\Delta \BYDEF \Delta(\epsilon)$ such that
\begin{equation*}
	\lim\limits_{\epsilon \rightarrow 0} \Delta(\epsilon) = 0 \quad \text{ and } \quad \lim\limits_{\epsilon \rightarrow 0} \frac{\Delta(\epsilon)}{\epsilon} = \infty
\end{equation*}
and set 
\begin{equation}
	\label{eq:tau_t}
	\tau_t \BYDEF t\Delta(\epsilon)  \quad t = 0, 1, \ldots, N(\epsilon) 
\end{equation}
where $N(\epsilon) = \left\lfloor \frac{1}{\Delta(\epsilon)}  \right\rfloor$.
\begin{Lemma}
	\label{lem:K_t}
	Define the function $K_t(\e) = \lfloor \tau_{t+1}/\e\rfloor-\lfloor \tau_{t}/\e\rfloor$. The following inequalities are valid:
	\begin{equation}
		\label{eq:first_K_t_inequality}
		\left|K_t(\e) - \frac{\Delta(\e)}{\e}\right| \le 1
	\end{equation}
	\begin{equation}
		\label{eq:second_K_t_inequality}
		\left|\frac{1}{K_t(\e)}- \frac{\e}{\Delta(\e)}\right| \le \left(\frac{\e}{\Delta(\e)}\right)^2\left(\frac{1}{1-\e/\Delta(\e)}\right)
	\end{equation}
	for all $t=0,1,\ldots, N(\e)-1$. 
\end{Lemma}
{\bf Proof.} The proof readily follows from the definition of the floor function (see also \cite{altman1997asymptotic} and \cite{brunetti2017nonzero}). $\ \Box$

\begin{Lemma}
	\label{lem:cts_to_discrete_fourth}
	There exists a constant $c > 0$ such that, for any $z \in Z$, any control plan $\pi \in \Pi$, and any sufficiently small $\e$,
	\begin{equation}
		\left|\left|\frac{1}{\Delta(\epsilon)}\int_{\tau_t}^{\tau_{t+1}}g(z,y^{\pi, y_0}_\e(s), u^{\pi, y_0}_\e(s))ds - \frac{1}{K(\epsilon)}\sum_{l=\lfloor \tau_t/\epsilon\rfloor}^{\lfloor \tau_{t}/\epsilon\rfloor+K(\epsilon)-1}g(z,y^{\pi, y_0}(l),u^{\pi, y_0}(l))\right|\right| \le c\left(\frac{\epsilon}{\Delta(\epsilon)}\right)
	\end{equation}
	where $(y^{\pi, y_0}(\cdot),u^{\pi, y_0}(\cdot))$ and $(y^{\pi, y_0}_{\epsilon}(\cdot),u^{\pi, y_0}_{\epsilon}(\cdot)) $ are defined in accordance with (\ref{eq:associated_system}) - (\ref{eq:control_dependency}) and (\ref{eq:hybrid_system_2}), and where
	\begin{equation}
		\label{eq:K_definition}
		K(\epsilon) \BYDEF \min\limits_{t=0,1,\ldots,N(\epsilon)-1}\left(\lfloor \tau_{t+1}/\epsilon\rfloor-\lfloor \tau_{t}/\epsilon\rfloor\right).
	\end{equation}
\end{Lemma}
{\bf Proof.} 
Let $\pi \in \Pi$ be a plan with the corresponding control $u^{\pi, y_0}(\cdot)$ and solution $y^{\pi, y_0}(\cdot)$. For any $z \in Z$, we have the following equality:
\begin{equation*}
	\begin{split}
		\frac{1}{\Delta(\epsilon)}\int_{\tau_t}^{\tau_{t+1}}g(z,y^{\pi, y_0}_\e(s), u^{\pi, y_0}_\e(s))ds = \ & \frac{1}{\Delta(\epsilon)}\int_{\e \lfloor \frac{\tau_t}{\e} \rfloor}^{\e \lfloor \frac{\tau_{t+1}}{\e} \rfloor}g(z,y^{\pi, y_0}_\e(s), u^{\pi, y_0}_\e(s))ds \\ &- \frac{1}{\Delta(\epsilon)}\int_{\e \lfloor \frac{\tau_t}{\e} \rfloor}^{\tau_{t}}g(z,y^{\pi, y_0}_\e(s), u^{\pi, y_0}_\e(s))ds\\ & + \frac{1}{\Delta(\epsilon)}\int_{\e \lfloor \frac{\tau_{t+1}}{\e} \rfloor}^{\tau_{t+1}}g(z,y^{\pi, y_0}_\e(s), u^{\pi, y_0}_\e(s))ds.
	\end{split}
\end{equation*}
By (\ref{eq:Lipschitz_bounded_assumption}),
\begin{equation*}
	\begin{split}
		\frac{1}{\Delta(\epsilon)}&\int_{\tau_t}^{\tau_{t+1}}g(z,y^{\pi, y_0}_\e(s), u^{\pi, y_0}_\e(s))ds \\\le & \ \frac{1}{\Delta(\epsilon)}\int_{\e \lfloor \frac{\tau_t}{\e} \rfloor}^{\e \lfloor \frac{\tau_{t+1}}{\e} \rfloor}g(z,y^{\pi, y_0}_\e(s), u^{\pi, y_0}_\e(s))ds + \frac{1}{\Delta(\epsilon)}\int_{\e \lfloor \frac{\tau_t}{\e} \rfloor}^{\tau_{t}}Mds + \frac{1}{\Delta(\epsilon)}\int_{\e \lfloor \frac{\tau_{t+1}}{\e} \rfloor}^{\tau_{t+1}}Mds \\
		&= \frac{1}{\Delta(\epsilon)}\int_{\e \lfloor \frac{\tau_t}{\e} \rfloor}^{\e \lfloor \frac{\tau_{t+1}}{\e} \rfloor}g(z,y^{\pi, y_0}_\e(s), u^{\pi, y_0}_\e(s))ds + \frac{M}{\Delta(\e)} \left(\tau_t - \e \left\lfloor \frac{\tau_t}{\e} \right\rfloor + \tau_{t+1} - \e \left\lfloor \frac{\tau_{t+1}}{\e}\right\rfloor\right)
		\\
		&\le \frac{1}{\Delta(\epsilon)}\int_{\e \lfloor \frac{\tau_t}{\e} \rfloor}^{\e \lfloor \frac{\tau_{t+1}}{\e} \rfloor}g(z,y^{\pi, y_0}_\e(s), u^{\pi, y_0}_\e(s))ds + 2M\left(\frac{\e}{\Delta(\e)}\right).
	\end{split}
\end{equation*}
Similarly, it can be shown that
\begin{equation*}
	\frac{1}{\Delta(\epsilon)}\int_{\e \lfloor \frac{\tau_t}{\e} \rfloor}^{\e \lfloor \frac{\tau_{t+1}}{\e} \rfloor}g(z,y^{\pi, y_0}_\e(s), u^{\pi, y_0}_\e(s))ds \le \frac{1}{\Delta(\epsilon)}\int_{\tau_t}^{\tau_{t+1}}g(z,y^{\pi, y_0}_\e(s), u^{\pi, y_0}_\e(s))ds + 2M\left(\frac{\e}{\Delta(\e)}\right),
\end{equation*}
implying 
\begin{equation}\label{eq:int_to_sum_eq_1}
	\left|\left| \frac{1}{\Delta(\epsilon)}\int_{\tau_t}^{\tau_{t+1}}g(z,y^{\pi, y_0}_\e(s), u^{\pi, y_0}_\e(s))ds - \frac{1}{\Delta(\epsilon)}\int_{\e \lfloor \frac{\tau_t}{\e} \rfloor}^{\e \lfloor \frac{\tau_{t+1}}{\e} \rfloor}g(z,y^{\pi, y_0}_\e(s), u^{\pi, y_0}_\e(s))ds \right|\right| \le 2M \left(\frac{\e }{\Delta(\e)}\right).
\end{equation}
Let us now consider the following equality:
\begin{equation*}
	\begin{aligned}
		\frac{1}{\Delta(\epsilon)}\int_{\e \lfloor \frac{\tau_t}{\e} \rfloor}^{\e \lfloor \frac{\tau_{t+1}}{\e} \rfloor}g(z,y^{\pi, y_0}_\e(s), u^{\pi, y_0}_\e(s))ds = \frac{1}{\Delta(\epsilon)}\int_{\e \lfloor \frac{\tau_t}{\e} \rfloor}^{\e( \lfloor \frac{\tau_{t}}{\e} \rfloor + K(\e))}g(z,y^{\pi, y_0}_\e(s), u^{\pi, y_0}_\e(s))ds   \\ 
		+  \ \frac{1}{\Delta(\epsilon)}\int_{\e (\lfloor \frac{\tau_t}{\e} \rfloor + K(\e))}^{\e\lfloor \frac{\tau_{t+1}}{\e} \rfloor}g(z,y^{\pi, y_0}_\e(s), u^{\pi, y_0}_\e(s))ds.
	\end{aligned}
\end{equation*}
By (\ref{eq:Lipschitz_bounded_assumption}), we have
\begin{align*}
		\frac{1}{\Delta(\epsilon)}\int_{\e \lfloor \frac{\tau_t}{\e} \rfloor}^{\e \lfloor \frac{\tau_{t+1}}{\e} \rfloor}g(z,y^{\pi, y_0}_\e(s), u^{\pi, y_0}_\e(s))ds &\le \frac{1}{\Delta(\epsilon)}\int_{\e \lfloor \frac{\tau_t}{\e} \rfloor}^{\e( \lfloor \frac{\tau_{t}}{\e} \rfloor + K(\e))}g(z,y^{\pi, y_0}_\e(s), u^{\pi, y_0}_\e(s))ds  \\
		&\ \ \ \ \ \ +\frac{1}{\Delta(\epsilon)}\int_{\e (\lfloor \frac{\tau_t}{\e} \rfloor + K(\e))}^{\e\lfloor \frac{\tau_{t+1}}{\e} \rfloor}Mds \\
		&= \frac{1}{\Delta(\epsilon)}\int_{\e \lfloor \frac{\tau_t}{\e} \rfloor}^{\e( \lfloor \frac{\tau_{t}}{\e} \rfloor + K(\e))}g(z,y^{\pi, y_0}_\e(s), u^{\pi, y_0}_\e(s))ds\\ 
		& \ \ \ \ \ \ + \frac{M}{\Delta(\epsilon)}\left(\e\left\lfloor \frac{\tau_{t+1}}{\e} \right\rfloor-\e \left(\left\lfloor \frac{\tau_t}{\e} \right\rfloor + K(\e)\right)\right)\\
		&\le \frac{1}{\Delta(\epsilon)}\int_{\e \lfloor \frac{\tau_t}{\e} \rfloor}^{\e( \lfloor \frac{\tau_{t}}{\e} \rfloor + K(\e))}g(z,y^{\pi, y_0}_\e(s), u^{\pi, y_0}_\e(s))ds\\
		& \ \ \ \ \ \ +\frac{M}{\Delta(\epsilon)}\left(\e\left(\frac{\tau_{t+1}}{\e}\right) -\e\left(\frac{\tau_t}{\e}\right) + \e - \e \left(\left(\frac{\tau_{t+1}}{\e}\right) - 1 - \left(\frac{\tau_t}{\e}\right)\right)\right)\\
		&= \frac{1}{\Delta(\epsilon)}\int_{\e \lfloor \frac{\tau_t}{\e} \rfloor}^{\e( \lfloor \frac{\tau_{t}}{\e} \rfloor + K(\e))}g(z,y^{\pi, y_0}_\e(s), u^{\pi, y_0}_\e(s))ds+2M\left(\frac{\e}{\Delta(\epsilon)}\right).
\end{align*}
Similarly,
\begin{equation*}
	\begin{split}
		\frac{1}{\Delta(\epsilon)}\int_{\e \lfloor \frac{\tau_t}{\e} \rfloor}^{\e( \lfloor \frac{\tau_{t}}{\e} \rfloor + K(\e))}g(z,y^{\pi, y_0}_\e(s), u^{\pi, y_0}_\e(s))ds &\le \frac{1}{\Delta(\epsilon)}\int_{\e \lfloor \frac{\tau_t}{\e} \rfloor}^{\e \lfloor \frac{\tau_{t+1}}{\e} \rfloor}g(z,y^{\pi, y_0}_\e(s), u^{\pi, y_0}_\e(s))ds \\
		& \ \ \  \ \ \ + 2M\left(\frac{\e}{\Delta(\epsilon)}\right),
	\end{split}
\end{equation*}
implying
\begin{equation}\label{eq:int_to_sum_eq_2}
	\begin{split}
		&\left|\left| \frac{1}{\Delta(\epsilon)}\int_{\e \lfloor \frac{\tau_t}{\e} \rfloor}^{\e \lfloor \frac{\tau_{t+1}}{\e} \rfloor}g(z,y^{\pi, y_0}_\e(s), u^{\pi, y_0}_\e(s))ds - \frac{1}{\Delta(\epsilon)}\int_{\e \lfloor \frac{\tau_t}{\e} \rfloor}^{\e( \lfloor \frac{\tau_{t}}{\e} \rfloor + K(\e))}g(z,y^{\pi, y_0}_\e(s), u^{\pi, y_0}_\e(s))ds \right|\right| \\ & \le 2M\left(\frac{\e}{\Delta(\epsilon)}\right).
	\end{split}
\end{equation}
Note that, by (\ref{eq:hybrid_system_2}),
\begin{equation}\label{eq:int_to_sum_eq_3}
	\frac{1}{\Delta(\epsilon)}\int_{\e \lfloor \frac{\tau_t}{\e} \rfloor}^{\e \lfloor \frac{\tau_{t}}{\e} \rfloor + K(\e)}g(z,y^{\pi, y_0}_\e(s), u^{\pi, y_0}_\e(s))ds = \frac{\e}{\Delta(\epsilon)}\sum_{l = \lfloor \frac{\tau_t}{\e}\rfloor }^{\lfloor \frac{\tau_t}{\e}\rfloor + K(\e)-1}g(z,y^{\pi, y_0}(l), u^{\pi, y_0}(l)).
\end{equation}
By (\ref{eq:Lipschitz_bounded_assumption}) and (\ref{eq:second_K_t_inequality}),
\begin{equation*}
	\begin{split}
		\left| \left| \frac{\e}{\Delta(\epsilon)}\sum_{l = \lfloor \frac{\tau_t}{\e}\rfloor }^{\lfloor \frac{\tau_t}{\e}\rfloor + K(\e)-1}g(z,y^{\pi, y_0}(l), u^{\pi, y_0}(l)) - \frac{1}{K(\e)}\sum_{l = \lfloor \frac{\tau_t}{\e}\rfloor }^{\lfloor \frac{\tau_t}{\e}\rfloor + K(\e)-1}g(z,y^{\pi, y_0}(l), u^{\pi, y_0}(l)) \right| \right| \\ \le \left(\frac{\e}{\Delta(\e)}\right)^2\left(\frac{1}{1-\frac{\e}{\Delta(\e)}}\right) K(\e) M.
	\end{split}
\end{equation*}
Since $\e$ was assumed to be small and $K(\e)$ is of the order $\Delta(\e) / \e$, we have 
\begin{equation}
	\label{eq:int_to_sum_eq_4}
	\begin{split}
		\left| \left| \frac{\e}{\Delta(\epsilon)}\sum_{l = \lfloor \frac{\tau_t}{\e}\rfloor }^{\lfloor \frac{\tau_t}{\e}\rfloor + K(\e)-1}g(z,y^{\pi, y_0}(l), u^{\pi, y_0}(l)) - \frac{1}{K(\e)}\sum_{l = \lfloor \frac{\tau_t}{\e}\rfloor }^{\lfloor \frac{\tau_t}{\e}\rfloor + K(\e)-1}g(z,y^{\pi, y_0}(l), u^{\pi, y_0}(l)) \right| \right| \\ \le c_1\left(\frac{\e}{\Delta(\e)}\right),
	\end{split}
\end{equation}
where $c_1 > 0$ is a constant. By combining (\ref{eq:int_to_sum_eq_1}), (\ref{eq:int_to_sum_eq_2}), (\ref{eq:int_to_sum_eq_3}) and  (\ref{eq:int_to_sum_eq_4}) and by setting $c = 4M+c_1$, one concludes the proof. $\ \Box$
\begin{Lemma}
	\label{lem:hybrid_system_soln_tightness}
	If $z_\e(\cdot)$ is a solution to the system (\ref{eq:hybrid_system}), then
	\begin{equation}
		||z_\e(s) - z_\e(\tau_t)|| \le M \Delta(\epsilon) \quad \forall \ s \in [\tau_t, \tau_{t+1}]
	\end{equation}
	for any range $[\tau_t, \tau_{t+1}]$ where $t=0,1,\ldots,N(\epsilon)-1$.
\end{Lemma}
{\bf Proof.} 
As a solution of the system (\ref{eq:hybrid_system}), $z_\e(s)$ allows the representation
\begin{equation*}
	z_\e(s) = z_\e(\tau_t) + \int_{\tau_t}^{s}g(z_\e(\tau), y_\e(\tau),u_\epsilon(\tau)) d\tau \quad \forall \ \tau \in [\tau_t, \tau_{t+1}].
\end{equation*}
Hence,
\begin{equation*}
	||z_\e(s) - z_\e(\tau_t)|| = \left|\left|\int_{\tau_t}^{s}g(z_\e(\tau),y_\e(\tau),u_\epsilon(\tau)) d\tau \right|\right| \quad \forall \ s \in [\tau_t, \tau_{t+1}],
\end{equation*}
thus, by (\ref{eq:Lipschitz_bounded_assumption}),
\begin{equation}
	||z_\e(s) - z_\e(\tau_t)|| \le \int_{\tau_t}^{s} Md\tau = M(s-\tau_t) \le M\Delta(\epsilon) \quad \forall \ s \in [\tau_t, \tau_{t+1}]
\end{equation}
as required. $\ \Box$
\begin{Lemma}
	\label{lem:Lipschitz_deterministic_system}
	If $\bar{z}(s)$ is a solution to the differential inclusion (\ref{eq:diff_inclusion}), then
	\begin{equation}
		\label{eq:zt_ztau_inequality}
		\left|\left|\bar{z}(s) - \bar{z}(\tau_t)\right|\right| \le M\Delta(\e) \quad \forall \ s \in [\tau_t, \tau_{t+1}]
	\end{equation}
	for any range $[\tau_t, \tau_{t+1}]$ where $t=0,1,\ldots,N(\e)-1$. 	
\end{Lemma}
{\bf Proof.} 
Let $\bar{z}(s)$ be a solution to the differential inclusion (\ref{eq:diff_inclusion}). By definition of $V_g(\bar{z}(s))$, there exists a $\g_s \in W$ such that
\begin{equation*}
	\dot{\bar{z}}(s) = \int_{Y \times \hat{U}} g(\bar{z}(s), y, u)\g_s(dy,du) \ \implies \ ||\dot{\bar{z}}(s) || \leq M,
\end{equation*}
the latter being implied by (\ref{eq:V_g_bounded}). Since for any $s \in [\tau_t, \tau_{t+1}]$,  $\bar{z}(s) = \bar{z}(\tau_t) + \int_{\tau_t}^{s}\dot{\bar{z}}(\tau)d\tau$, we have
\begin{equation*}
	\left|\left|\bar{z}(s) - \bar{z}(\tau_t)\right|\right| = \left|\left|\int_{\tau_t}^{s}\dot{\bar{z}}(\tau)d\tau\right|\right| \le \int_{\tau_t}^{s}\left|\left|\dot{\bar{z}}(\tau)\right|\right| d\tau
	\le \int_{\tau_t}^{s} M d\tau = M(s-\tau_t) \le M\Delta(\e),
\end{equation*}
proving the statement. $\ \Box$

\begin{Lemma}
	\label{lem:gronwall}
	Let $\ph(\epsilon) \ge 0$ be an arbitrary function of $\epsilon$ and let $M_\epsilon$ be a function of $\epsilon$ with values from the natural numbers, tending to infinity as $\epsilon$ tends to zero. Let $\Delta_0 = 0$ and $\Delta_t \ge 0$ satisfy the inequality $\Delta_{t+1} \le \Delta_t+L_1M_\epsilon^{-1}\Delta_t+\ph(\epsilon)M_\epsilon^{-1}$ for $t=1,\ldots,M_\epsilon - 1$, where $L_1$ is a positive constant. Then $\Delta_t$ also satisfies the inequality $\Delta_t \le \ph(\epsilon)L_1^{-1}e^{L_1}$ for $t=1,\ldots,M_\epsilon$.
\end{Lemma}
{\bf Proof.} See \cite[Proposition 5.1]{gaitsgory1992suboptimization} $\ \Box$\\
{\bf Proof of Theorem \ref{thm:first_hybrid_theorem}.}
Let $\pi \in \Pi$ be an arbitrary control plan and
let $(z^{\pi, y_0}_\e(s), y^{\pi, y_0}_\e(s), u^{\pi, y_0}_\e(s))$ be the corresponding solution and control of the system (\ref{eq:hybrid_system}). According to (\ref{eq:hybrid_system}), we have
\begin{equation}\label{eq:FHT-2}
	z^{\pi, y_0}_\e(\tau_{t+1}) = z^{\pi, y_0}_\e(\tau_t) +  \int_{\tau_t}^{\tau_{t+1}} g(z^{\pi, y_0}_\e(s),y^{\pi, y_0}_\e(s),u^{\pi, y_0}_\e(s))ds.
\end{equation}
Define a sequence of vectors $z_t$ by the equation
\begin{equation}\label{eq:FHT-1}
	z_{t+1} = z_t + \int_{\tau_t}^{\tau_{t+1}} g(z_t, y^{\pi, y_0}_\e(s),u^{\pi, y_0}_\e(s))ds
\end{equation}
for $t = 0, 1, \ldots, N(\epsilon) - 1$. 
By subtracting (\ref{eq:FHT-1}) from (\ref{eq:FHT-2}), we obtain
\begin{equation*}
	\begin{split}
		||z^{\pi, y_0}_\e(\tau_{t+1}) - z_{t+1}|| \le &||z^{\pi, y_0}_\e(\tau_t) - z_t|| \\ &+ \left| \left|\int_{\tau_t}^{\tau_{t+1}}g(z^{\pi, y_0}_\e(s),y^{\pi, y_0}_\e(s),u^{\pi, y_0}_\e(s))ds - \int_{\tau_t}^{\tau_{t+1}}g(z_t, y^{\pi, y_0}_\e(s),u^{\pi, y_0}_\e(s))ds\right|\right|.
	\end{split}
\end{equation*}
By the Lipschitz continuity of $g(\cdot)$ (see (\ref{eq:Lipschitz_g_assumption})), 
\begin{equation*}
	||z^{\pi, y_0}_\e(\tau_{t+1}) - z_{t+1}|| \le ||z^{\pi, y_0}_\e(\tau_t) - z_t|| + \int_{\tau_t}^{\tau_{t+1}}C(||z^{\pi, y_0}_\e(s) - z^{\pi, y_0}_\e(\tau_t)||+||z^{\pi, y_0}_\e(\tau_t) - z_t||)ds
\end{equation*}
which implies, by Lemma \ref{lem:hybrid_system_soln_tightness},
\begin{align*}
	||z^{\pi, y_0}_\e(\tau_{t+1}) - z_{t+1}|| & \le  ||z^{\pi, y_0}_\e(\tau_t) - z_t|| + \int_{\tau_t}^{\tau_{t+1}}(CM\Delta(\epsilon)+C||z^{\pi, y_0}_\e(\tau_t)-z_t||)ds \\
	& \le ||z^{\pi, y_0}_\e(\tau_t) - z_t||+ C\Delta(\epsilon)||z^{\pi, y_0}_\e(\tau_t) - z_t||+CM(\Delta(\epsilon))^2\\
	& \le ||z^{\pi, y_0}_\e(\tau_t) - z_t||+ C(1 / \left\lfloor 1 / \Delta(\e) \right\rfloor)||z^{\pi, y_0}_\e(\tau_t) - z_t||+CM\Delta(\epsilon)(1 / \left\lfloor 1 / \Delta(\e) \right\rfloor).
\end{align*}
Applying Lemma \ref{lem:gronwall} with $L_1 = C$, $M_\e = N(\e) = \left\lfloor 1 / \Delta(\e) \right\rfloor$ and $\ph(\epsilon) = CM\Delta(\epsilon)$ leads to the inequality
\begin{equation}
	\label{eq:Z_tau_t_z_t_inequality}
	||z^{\pi, y_0}_\e(\tau_t) - z_t|| \le M\Delta(\epsilon)e^{C}, \ \ t = 0,1,\ldots,N(\e).
\end{equation}
Define
\begin{equation}\label{eq:vl_definition}
	v_t = \text{argmin}\left\{ \left|\left|\frac{1}{K(\epsilon)}\sum_{l=\left\lfloor \frac{\tau_t}{\epsilon}\right\rfloor}^{\left\lfloor \frac{\tau_t}{\epsilon }\right\rfloor+K(\epsilon)-1}g(z_t,y(l),u(l)) - v \ \right|\right| \colon v \in V_g(z_t)  \right\}
\end{equation}
and the sequence of vectors
\begin{equation*}
	\zeta_{t+1} = \zeta_t + \Delta(\epsilon)\tilde{v}_t \quad t=0,1,\ldots,N(\epsilon) - 1,
\end{equation*}
where $\tilde{v}_t$ is the projection of $v_t$ onto the set $V_g(\zeta_t)$. By (\ref{eq:V_g_lipschitz_continuity}) and (\ref{eq:Z_tau_t_z_t_inequality}), we have
\begin{equation}
	\label{eq:v_t_tilde_v_t_difference}
	\begin{split}
		||v_t - \tilde{v}_t || &= d(v_t, V_g(\zeta_t)) \\ &\le d_H(V_g(z_t), V_g(\zeta_t)) \\ &\le  C||z_t - \zeta_t|| \\ &\le C||z_t - z^{\pi, y_0}_\e(\tau_t)|| \\ &+ C||z^{\pi, y_0}_\e(\tau_t) - \zeta_t|| \\ &\le CM\Delta(\e)e^{C}+C||z^{\pi, y_0}_\e(\tau_t) - \zeta_t||.
	\end{split}
\end{equation}
Consider now the following inequality:
\begin{equation}\label{eq:Ze_zeta_exp_difference}
	\begin{split}
		E[||z^{\pi, y_0}_\e(\tau_{t+1}) - \zeta_{t+1} ||] \le \ & E[||z^{\pi, y_0}_\e(\tau_t) - \zeta_t||] \\ & + E\left[\left|\left| \int_{\tau_t}^{\tau_{t+1}}g(z^{\pi, y_0}_\e(s), y^{\pi, y_0}_\e(s), u^{\pi, y_0}_\e(s), )ds - \Delta(\epsilon) \tilde{v}_t\right|\right|\right].
	\end{split}
\end{equation}
Note the second term on the right-hand side satisfies the inequality
\begin{equation}
	\label{eq:Ze_zeta_exp_difference_rhs}
	\begin{split}
		&E\left[\left|\left| \int_{\tau_t}^{\tau_{t+1}}g(z^{\pi, y_0}_\e(s), y^{\pi, y_0}_\e(s), u^{\pi, y_0}_\e(s))ds - \Delta(\e) \tilde{v}_t\right|\right|\right] \\ & \le E\left[  \int_{\tau_t}^{\tau_{t+1}}||g(z^{\pi, y_0}_\e(s), y^{\pi, y_0}_\e(s), u^{\pi, y_0}_\e(s)) - g(z^{\pi, y_0}_\e(\tau_t), y^{\pi, y_0}_\e(s), u^{\pi, y_0}_\e(s))|| ds\right.\\ & 
		\left. + \int_{\tau_t}^{\tau_{t+1}}||g(z^{\pi, y_0}_\e(\tau_t), y^{\pi, y_0}_\e(s), u^{\pi, y_0}_\e(s)) - g(z_t, y^{\pi, y_0}_\e(s), u^{\pi, y_0}_\e(s))||ds \right. \\ &  \left. + \left|\left|\int_{\tau_t}^{\tau_{t+1}}g(z_t, y^{\pi, y_0}_\e(s), u^{\pi, y_0}_\e(s))ds - \frac{ \Delta(\e)}{K(\e)}\sum_{l=\left\lfloor \frac{\tau_t}{\e}\right\rfloor+1}^{\left\lfloor \frac{\tau_t}{\e }\right\rfloor+K(\e)}g(z_t,y(l),u(l))  \right|\right| \right. \\ &
		\left. + \Delta(\e)\left|\left|\frac{ 1}{K(\e)}\sum_{l=\left\lfloor \frac{\tau_t}{\e}\right\rfloor}^{\left\lfloor \frac{\tau_t}{\e }\right\rfloor+K(\e)-1}g(z_t,y(l),u(l)) - v_t  \right|\right|+\Delta(\e)|| v_t - \ \tilde{v}_t ||\right].
	\end{split}
\end{equation}
Let us consider the right-hand side term-by-term. By (\ref{eq:Lipschitz_g_assumption}) and Lemma \ref{lem:hybrid_system_soln_tightness}, we have
\begin{equation}
	\label{eq:big_inequality_first}
	\begin{split}
		&E\left[\int_{\tau_t}^{\tau_{t+1}}||g(z^{\pi, y_0}_\e(s), y^{\pi, y_0}_\e(s), u^{\pi, y_0}_\e(s)) - g(z^{\pi, y_0}_\e(\tau_t), y^{\pi, y_0}_\e(s), u^{\pi, y_0}_\e(s))||ds\right]\\ &\le E\left[  C\int_{\tau_t}^{\tau_{t+1}}||z^{\pi, y_0}_\e(s)- z^{\pi, y_0}_\e(\tau_t)||ds\right] \\ 
		& \le  CM (\Delta(\e))^2
	\end{split}
\end{equation}
for the first term. By (\ref{eq:Lipschitz_g_assumption}) and  (\ref{eq:Z_tau_t_z_t_inequality}), we have
\begin{equation}
	\label{eq:big_inequality_second}
	\begin{split}
		&E\left[\int_{\tau_t}^{\tau_{t+1}}||g(z^{\pi, y_0}_\e(\tau_t), y^{\pi, y_0}_\e(s), u^{\pi, y_0}_\e(s)) - g(z_t, y^{\pi, y_0}_\e(s), u^{\pi, y_0}_\e(s))||ds\right]\\
		&\ \ \le E\left[C\int_{\tau_t}^{\tau_{t+1}}||z^{\pi, y_0}_\e(\tau_t)- z_t||ds\right] \\ 
		&\ \ \le \Delta(\e)CM\Delta(\e)e^{C}\\
		&\ \  = CM(\Delta(\e))^2e^{C}
	\end{split}
\end{equation}
for the second term. By Lemma \ref{lem:cts_to_discrete_fourth}, we have
\begin{equation}
	\label{eq:big_inequality_third}
	\left|\left|\int_{\tau_t}^{\tau_{t+1}}g(z_t, y^{\pi, y_0}_\e(s),u^{\pi, y_0}_\e(s))ds - \frac{ \Delta(\e)}{K(\e)}\sum_{m=\left\lfloor \frac{\tau_t}{\e}\right\rfloor}^{\left\lfloor \frac{\tau_t}{\e }\right\rfloor+K(\e)-1}g(z_t,y(m),u(m))  \right|\right| \le c\e
\end{equation}
for the third term. Now, since, 
\begin{equation*}
	\frac{1}{K(\e)}\sum_{m=\left\lfloor \frac{\tau_t}{\e}\right\rfloor}^{\left\lfloor \frac{\tau_t}{\e }\right\rfloor+K(\e)-1}g(z_t,y(m),u(m)) \in V_g\left(z_t, K(\e), y\left(\left\lfloor \frac{\tau_t}{\e} \right\rfloor\right)\right),
\end{equation*}
while $v_t$ is defined by (\ref{eq:vl_definition}), Corollary \ref{cor:vg_vgsy0_distance_expectation} implies for the fourth term  that
\begin{equation}
	\label{eq:big_inequality_fourth}
	\Delta(\e)E\left[\left|\left|\frac{ 1}{K(\e)}\sum_{m=\left\lfloor \frac{\tau_t}{\e}\right\rfloor}^{\left\lfloor \frac{\tau_t}{\e }\right\rfloor+K(\e)-1}g(z_t,y(m),u(m)) - v_t  \right|\right|\right] \le \Delta(\e)\bar{\nu}_g(K(\e)),
\end{equation}
where $\lim\limits_{\e \rightarrow 0}\bar{\nu}_g(K(\e)) = 0$. Finally, by (\ref{eq:v_t_tilde_v_t_difference}), we have 
\begin{equation}
	\label{eq:big_inequality_fifth}
	\Delta(\e)|| v_t - \ \tilde{v}_t || \le \Delta(\e)C||z^{\pi, y_0}_\e(\tau_t) - \zeta_t||+CM(\Delta(\e))^2e^{C}.
\end{equation}
Combining equations (\ref{eq:Ze_zeta_exp_difference_rhs}) - (\ref{eq:big_inequality_fifth}), we arrive at the inequality
\begin{equation}
	\label{eq:large_inequality}
	\begin{split}
		E&\left[\left|\left| \int_{\tau_t}^{\tau_{t+1}}g(z^{\pi_,y_0}_\e(s), y^{\pi_,y_0}_\e(s), u^{\pi_,y_0}_\e(s))ds - \Delta(\e) \tilde{v}_t\right|\right|\right] \\ & \le \Delta(\e)\big(CM \Delta(\e) + 2CM\Delta(\e)e^{C} + c\frac{\e}{\Delta(\e)} + \bar{\nu}_g(K(\e))  +  C||z^{\pi, y_0}_\e(\tau_t) - \zeta_t||\big).
	\end{split}
\end{equation}
Substituting back into inequality (\ref{eq:Ze_zeta_exp_difference}), we have
\begin{equation}
	\begin{split}
		&E[||z^{\pi_,y_0}_\e(\tau_{t+1}) - \zeta_{t+1} ||] \\ & \le ||z^{\pi, y_0}_\e(\tau_t)-\zeta_t||+C\Delta(\e)||z^{\pi, y_0}_\e(\tau_t)-\zeta_t|| + \Delta(\e)\left(CM \Delta(\e) + 2CM\Delta(\e)e^{C} + c\frac{\e}{\Delta(\e)}+ \bar{\nu}_g(K(\e))\right) \\
		& \le ||z^{\pi, y_0}_\e(\tau_t)-\zeta_t||+ C(1 / \lfloor 1 / \Delta(\e)\rfloor)||z^{\pi, y_0}_\e(\tau_t)-\zeta_t|| \\ & \ \ \ + (1 / \lfloor 1 / \Delta(\e)\rfloor)\left(CM \Delta(\e) + 2CMe^{C}\Delta(\e) + c\frac{\e}{\Delta(\e)}+ \bar{\nu}_g(K(\e))\right).
	\end{split}
\end{equation}
Applying Lemma \ref{lem:gronwall} with $L_1 = C$, $M_\e = \left\lfloor 1 / \Delta(\e) \right\rfloor$ and $\ph(\e) =  CM \Delta(\e) + 2CMe^{C}\Delta(\e) + c\e / \Delta(\e)+ \bar{\nu}_g(K(\e))$, we have
\begin{equation}
	\label{eq:Z_tau_t_and_zeta_t_difference}
	||z^{\pi, y_0}_\e(\tau_t)-\zeta_t|| \le \omega_1(\e).
\end{equation}
where the piece-wise linear function $\omega_1(\e) \BYDEF \frac{e^{C}}{C}(CM \Delta(\e) + 2CMe^{C}\Delta(\e) + c\e / \Delta(\e)+ \bar{\nu}_g(K(\e)))$ approaches zero as $\e \rightarrow 0$. Define now
\begin{equation}
	\label{eq:zeta(t)}
	\zeta(s) = \begin{cases}
		\zeta_t + (s-\tau_t)\tilde{v}_t \quad &\forall \ s \in [\tau_t, \tau_{t+1}], \ \ t \in \{0,1,\ldots,N(\e)-1\}\\
		\zeta_{N(\e)} + (s-\tau_{N(\e)})\tilde{v}_{N(\e)} \quad &\forall \ s \in [\tau_{N(\e)}, 1]
	\end{cases}
\end{equation}
and note that, for $s \in (\tau_t, \tau_{t+1})$,  $t=0,1,\ldots,N(\e)-1$:
\begin{align*}
	d(\dot{\zeta}(s), V_g(\zeta(s))) &= d(\tilde{v}_t, V_g(\zeta(s))) \\
	& \le d(\tilde{v}_t, V_g(\zeta_t)) + d_H(V_g(\zeta_t), V_g(\zeta(s))).
\end{align*}
Since $\tilde{v}_t \in V_g(\zeta_t)$, one has $d(\tilde{v}_t, V_g(\zeta_t)) = 0$. Also, by (\ref{eq:V_g_lipschitz_continuity}), 
\begin{align*}
	d_H(V_g(\zeta_t), V_g(\zeta(s))) &\le C||\zeta_t - \zeta(s) ||\\
	& = C(s-\tau_t) ||\tilde{v}_t|| \\
	& \le CM\Delta(\e)  
\end{align*}
for all $s \in [\tau_t, \tau_{t+1}]$,  $t=0,1,\ldots,N(\e)-1$. Similarly, 
\begin{equation*}
	d_H(V_g(\zeta_{N(\e)}), V_g(\zeta(s))) \le CM\Delta(\e) \quad \forall \ s \in [\tau_{N(\e)}, 1].
\end{equation*}
Hence,
\begin{equation}
	d(\dot{\zeta}(s), V_g(\zeta(s))) \le CM\Delta(\e) \quad \text{for almost all} \ s \in [0, 1].
\end{equation}
By Filippov's Theorem (see \cite[p. ~171]{aubin2009viability})
, there exists $z_\e(s)$ solving the differential inclusion (\ref{eq:diff_inclusion}) such that
\begin{equation}
	\label{eq:filippov_solution}
	\max_{s \in [0, 1]}||z_\e(s) - \zeta(s)|| \le C_1 \Delta(\e).
\end{equation}
for some $C_1 > 0$. By (\ref{eq:filippov_solution}) and (\ref{eq:Z_tau_t_and_zeta_t_difference}),
\begin{equation}
	\label{eq:za_tt_Zsp_tau_t_difference_equation}
	\begin{split}
		E[||z_\e(\tau_t) - z^{\pi, y_0}_\e(\tau_t)||] &\le E[||z_\e(\tau_t) - \zeta(\tau_t)||] + E[||\zeta(\tau_t) - \zeta_t||] + E[||\zeta_t - z^{\pi, y_0}_\e(\tau_t)||] \\
		& \le C_1 \Delta(\e) + M\Delta(\e) + \omega_1(\e),
	\end{split}
\end{equation}
where $\omega_1(\e)$ is as defined in (\ref{eq:Z_tau_t_and_zeta_t_difference}). Finally, by Lemma \ref{lem:hybrid_system_soln_tightness}, (\ref{eq:za_tt_Zsp_tau_t_difference_equation}) and Lemma \ref{lem:Lipschitz_deterministic_system},
\begin{align*}
	E[||z^{\pi, y_0}_\e(s) - z_\e(s)||] &\le E[||z^{\pi, y_0}_\e(s) - z^{\pi, y_0}_\e(\tau_t)||] + E[||z^{\pi, y_0}_\e(\tau_t) - z_\e(\tau_t)||] + E[||z_\e(\tau_t) - z_\e(s)||] \\
	&\le C_1\Delta(\e) + 3M\Delta(\e) + \omega_1(\e),
\end{align*}
which completes the proof. $\ \Box$



\section{Construction of control plans that generate solutions of the hybrid system, the $z$-components of which are close to the solutions of the differential inclusion}
\label{sec:hybrid_system_convergence}
For brevity, throughout this section we do not emphasise the dependence of state-control trajectories on initial conditions.

\begin{Theorem}
	\label{thm:second_hybrid_theorem}
	Assume the system (\ref{eq:associated_system}) satisfies the strong $h$-approximation conditions (see Definition \ref{def:S-h-AC}) for any vector function $h(\cdot,\cdot)$ as in (\ref{eq:h-AC_function}). Then, corresponding to any solution $\bar{z}(\cdot)$ of the differential inclusion (\ref{eq:diff_inclusion}), there exists a control plan $\bar{\pi}$ such that the solution $z^{\bar{\pi}}_{\epsilon}(\cdot)$ of system (\ref{eq:hybrid_system}) obtained with the use of this control plan satisfies the inequality
	\begin{equation}\label{eq:second_hybrid_theorem_eqn}
		\max_{s \in [0, 1]} E ||z^{\bar{\pi}}_{\epsilon}(s) - \bar{z}(s)|| \le \psi_2(\epsilon)
	\end{equation}
	for some $\psi_2(\epsilon)$ such that $\lim\limits_{\e \rightarrow 0}\psi_2(\epsilon) = 0$. 	
\end{Theorem}
{\bf Proof.} The proof of this theorem is given at the end of this section. $\ \Box$
\begin{Corollary}\label{cor:hybrid_main_corollary}
	Assume the system (\ref{eq:associated_system}) satisfies the strong $h$-approximation conditions (see Definition \ref{def:S-h-AC}) for any vector function $h(\cdot,\cdot)$ as in (\ref{eq:h-AC_function}). Then, the optimal value $F^\epsilon$ of problem (\ref{eq:OC_hybrid}) converges to the optimal value $F^0$ of the problem (\ref{eq:OC_diff_inclusion}) as $\epsilon \rightarrow 0$. More precisely,
	\begin{equation}\label{eq:hybrid_corollary_1}
		|F^\epsilon - F^0 | \le C_G\psi(\epsilon),
	\end{equation}
	where $\psi(\epsilon) =\max(\psi_1(\epsilon),\psi_2(\epsilon))$ (with $\psi_1(\epsilon)$ being as in (\ref{eq:DI_sol_existence}) and $\psi_2(\epsilon)$ being as in (\ref{eq:second_hybrid_theorem_eqn})). 
	In addition, if $\bar{z}(s)$ is an optimal solution for the problem (\ref{eq:OC_diff_inclusion}), the control plan $\bar{\pi}\in \Pi$ constructed in the proof of Theorem \ref{thm:second_hybrid_theorem} generates the solution $z^{\bar{\pi}}_\e(\cdot)$ of system (\ref{eq:hybrid_system_1}) such that
	\begin{equation}\label{eq:hybrid_corollary_2}
		E[G(z^{\bar{\pi}}_\e(1))]\leq F^\epsilon + C_G(\psi(\e) + \psi_2(\e)).
	\end{equation}
	That is, the control plan $\bar{\pi}$ is asymptotically optimal for the problem (\ref{eq:OC_hybrid}).
\end{Corollary}
{\bf Proof.} Let $\pi \in \Pi$ be an arbitrary control plan, let $z^{\pi}_{\e}(\cdot)$ be the corresponding solution of (\ref{eq:hybrid_system})
and let $z_\e(\cdot)$ be the solution to the differential inclusion (\ref{eq:diff_inclusion}) that satisfies (\ref{eq:DI_sol_existence}). Then,
\begin{equation}
	\label{eq:main_corollary_eq_1}
	|E[G(z^{\pi}_\e(1))] - E[G(z_\e(1))]| \le E[C_G|z^{\pi}_\e(1) - z_\e(1)|] \le C_G \psi_1(\epsilon),
\end{equation}
where $C_G$ is the Lipschitz constant of $G(\cdot)$. Note that
\begin{equation*}
	G(z_\e(1)) \ge F^0 \implies
	E[G(z_\e(1))]  = G(z_\e(1)) \ge F^0,
\end{equation*}
implying, by (\ref{eq:main_corollary_eq_1}),
\begin{equation*}
	E[G(z^{\pi}_\e(1))] \ge E[G(z_\e(1))] - C_G\psi_1(\e) \ge F^0 - C_G\psi_1(\e).
\end{equation*}
Since $\pi$ was an arbitrarily selected control plan, 
\begin{equation}
	\label{eq:main_corollary_eq_2}
	F^\epsilon = \inf_{\pi \in \Pi}E[G(z^{\pi}_\e(1))] \ge F^0 - C_G\psi_1(\epsilon).
\end{equation}
Now, let $\bar{z}(\cdot)$ be an optimal trajectory	for the problem (\ref{eq:OC_diff_inclusion}). By Theorem \ref{thm:second_hybrid_theorem}, there exists a control plan $\bar{\pi}$ with the corresponding solution $z^{\bar{\pi}}_{\epsilon}(\cdot)$ obtained from the system (\ref{eq:hybrid_system}) such that
\begin{equation}\label{eq:main_corollary_eq_3}
	\begin{split}
		|E[G(z^{\bar{\pi}}_{\epsilon}(1))] - F^0 | &= |E[G(z^{\bar{\pi}}_{\epsilon}(1))] - E[G(\bar{z}(1))]| \\
		&\le C_GE[||z^{\bar{\pi}}_{\epsilon}(1) - \bar{z}(1)||]\\ 
		&\le C_G \psi_2(\epsilon),
	\end{split}
\end{equation}
implying
\begin{equation}
	\label{eq:main_corollary_eq_4}
	E[G(z^{\bar{\pi}}_{\epsilon}(1))] \le F^0 + C_G\psi_2(\epsilon).
\end{equation}
Since $E[G(z^{\bar{\pi}}_{\epsilon}(1))] \ge F^\epsilon$, we have by (\ref{eq:main_corollary_eq_4}) that
\begin{equation}\label{eq:main_corollary_eq_5}
	F^{\epsilon}\leq F^0 +C_G\psi_2(\epsilon).
\end{equation}
This along with (\ref{eq:main_corollary_eq_2}) proves (\ref{eq:hybrid_corollary_1}). The validity of (\ref{eq:hybrid_corollary_2}) follows from  (\ref{eq:hybrid_corollary_1}) and (\ref{eq:main_corollary_eq_3}).  $\ \Box$\\
Let us construct the control plan $\bar{\pi}$ providing the inequality (\ref{eq:second_hybrid_theorem_eqn}). Let $\bar{z}(\cdot)$ be a solution to the differential inclusion (\ref{eq:diff_inclusion}) and define the projection of $\left(\Delta(\e)^{-1} \int_{0}^{\tau_{1}}\dot{\bar{z}}(s)ds \right)$ onto the set $V_g(\bar{z}(0))$ as follows:
\begin{equation}
	\label{eq:projection_vector_vz_0}
	v_0 \BYDEF \text{argmin}\left\{\left|\left|\Delta(\e)^{-1} \int_{0}^{\tau_{1}}\dot{\bar{z}}(s) ds - v\right|\right| : v \in V_g(\bar{z}(0)) \right\}.	
\end{equation}
Then, by Corollary \ref{cor:vg_vgsy0_distance}, for any initial condition $y_0 \in Y$, there exists a control plan $$\pi^0=\big(\pi^0_0, \pi^0_1(\sigma(0)), \ldots, \pi^0_{K(\e)-1}(\sigma(0), \ldots, \sigma(K(\e)-2))\big)$$  with the corresponding state-control trajectory $(y^{\pi^0}(\cdot), u^{\pi^0}(\cdot))$ that generates
the random occupational measure $\chi^{\pi^0, y_0, K(\e)} \in B_{K(\e)}(y_0)$ such that
\begin{equation}
	E\left[\left|\left|v_0 -  \int_{Y \times \hat{U}}^{} g(\bar{z}(0), y, u) \chi^{\pi^0, y_0, K(\e)}(dy, du)\right|\right|\right] \le \bar{\nu}_g(K(\e)),
\end{equation}
where $\lim\limits_{\e \rightarrow 0}\bar{\nu}_g\left(K(\e)\right) = 0$, or equivalently (see (\ref{G88})),
\begin{equation*}
	E\left[\left|\left|v_0 -  \frac{1}{K(\e)}\sum_{m = 0}^{K(\e)-1} g(\bar{z}(0), y^{\pi^0}(m), u^{\pi^0}(m)) \right|\right|\right] \le \bar{\nu}_g(K(\e)).
\end{equation*}
We begin the construction of the control plan $\bar{\pi} = (\bar{\pi}_k(\cdot), \ k = 0,1,\ldots, \lfloor 1 / \e \rfloor) \in \Pi$ by defining it for $k = 0,1,\ldots,\lfloor \tau_1 / \e \rfloor -1$ as follows:
$$\bar \pi_k(\cdot) = \pi^0_k(\cdot) \ \text{for} \  k = 0,1,\ldots,K(\e)-1, \ \text{and} \  \bar \pi_k(\cdot) =u \ \text{for} \
k=K(\e),\ldots, [\tau_1 / \e]-1,$$
where $u$ is an arbitrary element of $\hat U$. Let us assume that the control plan $\bar{\pi}$ has been constructed for the time moments $k = 0,1,\ldots,\lfloor \tau_t / \e \rfloor -1$ and, by use of induction, extend the construction to the sequence of time moments $k = \lfloor \tau_t / \e \rfloor, \ldots, \lfloor \tau_{t+1} / \e \rfloor - 1$. 

Define $v_t$ as the projection of the vector $\left(\Delta(\e)^{-1} \int_{\tau_{t}}^{\tau_{t+1}}\dot{\bar{z}}(s)ds \right)$ onto the set $V_g(\bar{z}(\tau_t))$:
\begin{equation}
	\label{eq:projection_vector_vz_t}
	v_t \BYDEF \text{argmin}\left\{\left|\left|\Delta(\e)^{-1} \int_{\tau_t}^{\tau_{t+1}}\dot{\bar{z}}(s) ds - v\right|\right| : v \in V_g(\bar{z}(\tau_t)) \right\}.
\end{equation}
By Corollary \ref{cor:vg_vgsy0_distance}, for the initial condition $y(\lfloor \tau_t / \e \rfloor) = y^t$, there exists a control plan
$$\pi^t=\big(\pi^t_{\lfloor \tau_t / \e \rfloor}, \pi^t_{\lfloor \tau_t / \e \rfloor+1}(\sigma(\lfloor \tau_t / \e \rfloor)),\ldots, \pi^t_{\lfloor \tau_{t+1} / \e \rfloor-1}(\sigma(\lfloor \tau_t / \e \rfloor), \ldots, \sigma(\lfloor \tau_{t} / \e \rfloor+K(\e)-2))\big)$$
with the corresponding state-control trajectory $(y^{\pi^t}(\cdot), u^{\pi^t}(\cdot))$ that generates
the occupational measure $\chi^{\pi^t, y^t, K(\e)} \in B_{K(\e)}(\lfloor \tau_t / \e \rfloor)$ such that
\begin{equation}
	E\left[\left|\left|v_t -  \int_{Y \times \hat{U}}^{} g(\bar{z}(\tau_t), y, u) \chi^{\pi^t, y^t, K(\e)}(dy, du)\right|\right|\right] \le \bar{\nu}_g(K(\e)),
\end{equation}
where $\lim\limits_{\e \rightarrow 0}\bar{\nu}_g\left(K(\e)\right) = 0$, or equivalently,
\begin{equation}\label{eq:v_t_measures_expected_difference}
	E\left[\left|\left|v_t -  \frac{1}{K(\e)}\sum_{m = \lfloor \tau_t / \e \rfloor}^{\lfloor \tau_t / \e \rfloor + K(\e)-1} g(\bar{z}(\tau_t), y^{\pi^t}(m), u^{\pi^t}(m)) \right|\right|\right] \le \bar{\nu}_g(K(\e)).
\end{equation}
Extend the definition of the control plan $\bar{\pi}$ to the sequence $k = \lfloor\tau_t / \e \rfloor, \ldots,  \lfloor\tau_{t+1} / \e \rfloor-1$ as follows:
\begin{equation}\label{eq:corr_cp_hybrid_1}
	\begin{split}
		\bar{\pi}_k(\cdot) = \pi_k^t(\cdot) \ \text{for} \ k=\lfloor\tau_t / \e\rfloor, \lfloor \tau_t / \e\rfloor+1, \ldots, \lfloor \tau_t / \e\rfloor+K(\e)-1,\\
		\bar{\pi}_k(\cdot) = u \ \text{for} \ k = \lfloor \tau_t / \e\rfloor+K(\e), \ldots, \lfloor \tau_{t+1} / \e\rfloor-1,
	\end{split}
\end{equation}
where $u$ is an arbitrary element of $\hat U$. By repetitions of the preceding process, we define $\bar{\pi}$ for $k=0,1,\ldots, \lfloor  \tau_{N(\e)} / \e \rfloor-1$. We complete the definition of the control plan $\bar{\pi}$ by defining it for $k=\lfloor  \tau_{N(\e)} / \e \rfloor, \ldots, \lfloor 1 / \e \rfloor$ by taking
\begin{equation}
	\label{eq:corr_cp_hybrid_2}
	\bar{\pi}_k(\cdot) = u \ \text{for} \ k=\lfloor  \tau_{N(\e)} / \e \rfloor, \ldots, \lfloor 1 / \e \rfloor,
\end{equation}
where $u$ is an arbitrary element of $\hat U$. 

The fact that the control plan $\bar{\pi}$ provides inequality (\ref{eq:second_hybrid_theorem_eqn}) is based on the following lemmas.

\begin{Lemma}
	\label{lem:V_distance}
	Let $\bar{z}(s)$ be a solution to the differential inclusion (\ref{eq:diff_inclusion}). Then
	\begin{equation*}
		d\left(\Delta(\e)^{-1}\int_{\tau_t}^{\tau_{t+1}}{\dot{\bar{z}}}(s)ds, V_g(\bar{z}(\tau_t))\right) \le CM\Delta(\e)
	\end{equation*}
	where the distance $d(\cdot, \cdot)$ is as introduced in (\ref{def:distance}) and $C$ and $M$ are the positive constants introduced with (\ref{eq:Lipschitz_g_assumption}) and (\ref{eq:Lipschitz_bounded_assumption}), respectively.
	
\end{Lemma}
{\bf Proof.} 
Let $\bar{z}(s)$ be a solution to the differential inclusion (\ref{eq:diff_inclusion}). On the interval $s \in [\tau_t, \tau_{t+1}]$ we have, by (\ref{eq:V_g_lipschitz_continuity}) and Lemma \ref{lem:Lipschitz_deterministic_system}, 
\begin{align*}
	\dot{\bar{z}}(s)\in V_g(\bar{z}(s)) &\subset V_g(\bar{z}(\tau_t))+C||\bar{z}(s) - \bar{z}(\tau_t)||\bar{B}
	\\ &\subset V_g(\bar{z}(\tau_t))+CM\Delta(\e)\bar{B}
\end{align*}   
for the closed unit ball $\bar{B} \in \mathbb{R}^k$. By the convexity of $V_g(\bar{z})$, 
\begin{equation*}
	\left(\Delta(\e)^{-1}\int_{\tau_t}^{\tau_{t+1}}{\dot{\bar{z}}}(s)ds\right) \in V_g(\bar{z}(\tau_t))+CM\Delta(\e)\bar{B}
\end{equation*}
implying
\begin{equation*}
	d\left(\Delta(\e)^{-1}\int_{\tau_t}^{\tau_{t+1}}{\dot{\bar{z}}}(s)ds, V_g(\bar{z}(\tau_t))\right) \le CM\Delta(\e)
\end{equation*}
and completing the proof. $\ \Box$

\begin{Lemma}
	Let $\bar{z}(s)$ be an arbitrary solution to the differential inclusion (\ref{eq:diff_inclusion}) and let $\zeta_t$,  $t=0,1,\ldots,N(\e) -1$ be defined by the equations:
	\begin{equation}
		\label{eq:zeta_t_difference}
		\zeta_{t+1} = \zeta_t + \Delta(\e)v_t, \quad t=0,1,\ldots,N(\e) - 1, \quad \zeta_0=z_0,
	\end{equation}
	where the vector $v_t$ is as defined in (\ref{eq:projection_vector_vz_0}) for $t=0$, and in (\ref{eq:projection_vector_vz_t}) for $t>0$. Then
	\begin{equation}
		\label{eq:bar_z_tau_t_and_zeta_t}
		\left|\left|\bar{z}(\tau_{t})-\zeta_{t}\right|\right| \le CM\Delta(\e)
	\end{equation}
	for all $t = 0,1,\ldots, N(\e).$
\end{Lemma}
{\bf Proof.} 
Let $\bar{z}(s)$ be an arbitrary solution to the differential inclusion (\ref{eq:diff_inclusion}). By Lemma \ref{lem:V_distance},
\begin{align*}
	\left|\left|\bar{z}(\tau_{t+1})-\zeta_{t+1}\right|\right| &\le \left|\left|\bar{z}(\tau_{t})-\zeta_{t}\right|\right| + \Delta(\e)\left|\left|\Delta^{-1}(\e)\int_{\tau_t}^{\tau_{t+1}}\dot{\bar{z}}(s)ds - v_t\right|\right| \\ 
	& \le \left|\left|\bar{z}(\tau_{t})-\zeta_{t}\right|\right| + CM\Delta(\e)^2,
\end{align*}
for all $t = 0, 1, \ldots, N(\e) - 1$. It can be readily seen that
\begin{equation}
	\left|\left|\bar{z}(\tau_{t})-\zeta_{t}\right|\right| \le N(\e)CM(\Delta(\e))^2 \le CM\Delta(\e), \ \ \ t = 0,1,\ldots, N(\e), 
\end{equation}
as required. $\ \Box$ \\
{\bf Proof of Theorem \ref{thm:second_hybrid_theorem}.} Let $\bar{\pi}$ be the control plan from (\ref{eq:corr_cp_hybrid_1}) - (\ref{eq:corr_cp_hybrid_2}) and define the sequence of random variables $z_t$ for $t=0,1,\ldots,N(\e)$ according to the following:
\begin{equation}
	\label{eq:z_t_difference}
	z_{t+1} = z_t + \int_{\tau_t}^{\tau_{t+1}}g(z_t, y^{\bar{\pi}}_\e(s), u^{\bar{\pi}}_\e(s))ds, \quad Z_0 = z_0,
\end{equation}
where $u^{\bar{\pi}}_\e(s)$ and $y^{\bar{\pi}}_\e(s)$ are the controls and states corresponding to the control plan $\bar{\pi}$, obtained from (\ref{eq:hybrid_system_2}). By taking the difference between (\ref{eq:z_t_difference}) and (\ref{eq:zeta_t_difference}), we obtain
\begin{equation}
	\label{eq:zeta_Z_difference_equation}
	\begin{split}
		E&\left[\left|\left|\zeta_{t+1} - z_{t+1}\right|\right|\right] \le  E\left[\left|\left|\zeta_{t} - z_{t}\right|\right|\right] + \Delta(\e)E\left[\left|\left|\frac{1}{\Delta(\e)}\int_{\tau_t}^{\tau_{t+1}}g(z_t, y^{\bar{\pi}}_\e(s), u^{\bar{\pi}}_\e(s))ds - v_t\right|\right|\right] \\&
		\le E\left[\left|\left|\zeta_{t} - z_{t}\right|\right|\right] + \Delta(\e)E\left[\left|\left|\frac{1}{\Delta(\e)}\int_{\tau_t}^{\tau_{t+1}}\left(g(z_t, y^{\bar{\pi}}_\e(s), u^{\bar{\pi}}_\e(s)) - g(\bar{z}(\tau_t), y^{\bar{\pi}}_\e(s),u^{\bar{\pi}}_\e(s))\right)ds\right|\right|\right] \\ & \ \ \ \ + \Delta(\e)E\left[\left|\left|\frac{1}{\Delta(\e)}\int_{\tau_t}^{\tau_{t+1}}g(\bar{z}(\tau_t), y^{\bar{\pi}}_\e(s), u^{\bar{\pi}}_\e(s))ds - v_t\right|\right|\right].
	\end{split}
\end{equation}
By (\ref{eq:Lipschitz_g_assumption}) and (\ref{eq:bar_z_tau_t_and_zeta_t}),
\begin{equation}
	\label{eq:zeta_Z_difference_equation_first_result}
	\begin{split}
		E&\left[\left|\left|\frac{1}{\Delta(\e)}\int_{\tau_t}^{\tau_{t+1}}\left(g(z_t, y^{\bar{\pi}}_\e(s), u^{\bar{\pi}}_\e(s)) - g(\bar{z}(\tau_t), y^{\bar{\pi}}_\e(s),u^{\bar{\pi}}_\e(s)\right)ds\right|\right|\right] \\ &\le E\left[\frac{1}{\Delta(\e)}\int_{\tau_t}^{\tau_{t+1}}C\left|\left|z_t - \bar{z}(\tau_{t})\right|\right|ds\right] \\
		& \le E\left[\frac{1}{\Delta(\e)}\int_{\tau_t}^{\tau_{t+1}}C(\left|\left|\zeta_t - z_t\right|\right| + \left|\left|\zeta_t - \bar{z}(\tau_{t})\right|\right|) ds\right] \\
		& \le E\left[C\left|\left|\zeta_t - z_t\right|\right| + C\left|\left|\zeta_t - \bar{z}(\tau_{t})\right|\right| \right] \\
		& \le CE\left[\left|\left|\zeta_t - z_t\right|\right|\right] + C^2TM\Delta(\e).  
	\end{split}
\end{equation}
By Lemma \ref{lem:cts_to_discrete_fourth},
\begin{align*}
	&\left|\left|\frac{1}{\Delta(\e)}\int_{\tau_t}^{\tau_{t+1}}g(\bar{z}(\tau_t), y^{\bar{\pi}}_\e(s), u^{\bar{\pi}}_\e(s))ds - v_t\right|\right| \\
	&\le \left|\left|\frac{1}{\Delta(\e)}\int_{\tau_t}^{\tau_{t+1}}g(\bar{z}(\tau_t), y^{\bar{\pi}}_\e(s), u^{\bar{\pi}}_\e(s))ds - \frac{1}{K(\e)}\sum_{m=\lfloor \tau_t/\e\rfloor}^{\lfloor \tau_{t}/\e\rfloor+K(\e)-1}g(\bar{z}(\tau_t), y^{\bar{\pi}}_\e(m),u^{\bar{\pi}}_\e(m))\right|\right|\\ & +\left|\left|\frac{1}{K(\e)}\sum_{m=\lfloor \tau_t/\e\rfloor}^{\lfloor \tau_{t}/\e\rfloor+K(\e)-1}g(\bar{z}(\tau_t), y^{\bar{\pi}}_\e(m),u^{\bar{\pi}}_\e(m)) - v_t\right|\right| \\
	& \le c\left(\frac{\e}{\Delta(\e)}\right) + \left|\left|\frac{1}{K(\e)}\sum_{m=\lfloor \tau_t/\e\rfloor}^{\lfloor \tau_{t}/\e\rfloor+K(\e)-1}g(\bar{z}(\tau_t), y^{\bar{\pi}}(m), u^{\bar{\pi}}(m)) - v_t\right|\right|
\end{align*}
for some constant $c > 0$. Hence, by the construction of $\bar{\pi}$ (in particular, see (\ref{eq:v_t_measures_expected_difference})),
\begin{equation}
	\label{eq:zeta_Z_difference_equation_second_result}
	\left|\left|\frac{1}{\Delta(\e)}\int_{\tau_t}^{\tau_{t+1}}g(\bar{z}(\tau_t), y^{\bar{\pi}}_\e(s), u^{\bar{\pi}}_\e(s))ds - v_t\right|\right| \le c\left(\frac{\e}{\Delta(\e)}\right) + \bar{\nu}_g(K(\e)).
\end{equation}
By (\ref{eq:zeta_Z_difference_equation_first_result}) and (\ref{eq:zeta_Z_difference_equation_second_result}), (\ref{eq:zeta_Z_difference_equation}) becomes
\begin{equation*}
	\begin{split}
		E\left[\left|\left|\zeta_{t+1} - z_{t+1}\right|\right|\right] &\le  E\left[\left|\left|\zeta_{t} - z_{t}\right|\right|\right] + \Delta(\e)\left(CE\left[\left|\left|\zeta_t - z_t\right|\right|\right] + C^2M\Delta(\e) + c\left(\frac{\e}{\Delta(\e)}\right) + \bar{\nu}_g(K(\e))\right).\\&
		\le E\left[\left|\left|\zeta_{t} - z_{t}\right|\right|\right] + C(1/\lfloor 1 / \Delta(\e) \rfloor)E\left[\left|\left|\zeta_t - z_t\right|\right|\right] \\
		& \ \ \ + (1/\lfloor 1 / \Delta(\e) \rfloor)\left(C^2M\Delta(\e) + c\left(\frac{\e}{\Delta(\e)}\right) + \bar{\nu}_g(K(\e))\right).
	\end{split}
\end{equation*}
Applying Lemma \ref{lem:gronwall}, with $L_1 = C$, $M_\e = \lfloor 1/ \Delta(\e) \rfloor$ and $\ph(\e) = \left(C^2M\Delta(\e) + c\left(\frac{\e}{\Delta(\e)}\right)+\bar{\nu}_g(K(\e))\right)$, we obtain
\begin{equation}
	\label{eq:zeta_t_and_z_t}
	E\left[\left|\left|\zeta_{t} - z_{t}\right|\right|\right] \le \frac{1}{C}\left(C^2M\Delta(\e) +c\left(\frac{\e}{\Delta(\e)}\right) +  \bar{\nu}_g(K(\e))\right)e^{C}, \ \ \  t=0,1,\ldots,N(\e).
\end{equation}
Now, since $z^{\bar{\pi}}_\e(\cdot)$ is the solution of (\ref{eq:hybrid_system}) associated with the control $u^{\bar{\pi}}_\e(\cdot)$,
\begin{equation*}
	z^{\bar{\pi}}_\e(\tau_{t+1}) = z^{\bar{\pi}}_\e(\tau_t) + \int_{\tau_{t}}^{\tau_{t+1}} g(z^{\bar{\pi}}_\e(s),y^{\bar{\pi}}_{\e}(s),u^{\bar{\pi}}_{\e}(s))ds.
\end{equation*}
Consider the difference between the above expression and (\ref{eq:z_t_difference}):
\begin{equation}
	\label{eq:difference_equation_Z_tau_t+1_z_t+1}
	\begin{split}
		E\left[\left|\left|z^{\bar{\pi}}_\e(\tau_{t+1})-z_{t+1}\right|\right|\right] &\le E\left[||z^{\bar{\pi}}_\e(\tau_t)-z_t|| + \int_{\tau_{t}}^{\tau_{t+1}} C||z^{\bar{\pi}}_\e(s)-z_t||ds\right] \\
		&\le E\left[||z^{\bar{\pi}}_\e(\tau_t)-z_t|| + \int_{\tau_{t}}^{\tau_{t+1}} (C||z^{\bar{\pi}}_\e(s)-z^{\bar{\pi}}_\e(\tau_t)|| + C||z^{\bar{\pi}}_\e(\tau_t)-z_t||)ds\right].
	\end{split} 
\end{equation}
By Lemma \ref{lem:hybrid_system_soln_tightness},
\begin{align*}
	E&[||z^{\bar{\pi}}_\e(\tau_{t+1})-z_{t+1}||] \le E[(1+C\Delta(\e))||z^{\bar{\pi}}_\e(\tau_t)-z_t|| + CM(\Delta(\e))^2] \\
	&\le E[||z^{\bar{\pi}}_\e(\tau_t)-z_t||] + C(1/\lfloor 1/\Delta(\e)\rfloor)E[||z^{\bar{\pi}}_\e(\tau_t)-z_t||] +  CM\Delta(\e)(1/\lfloor 1/\Delta(\e)\rfloor)
\end{align*}
and by Lemma \ref{lem:gronwall} with $L_1 = C$, $M_\e = \lfloor 1/\Delta(\e)\rfloor$ and $\ph(\e) = CM\Delta(\e)$,
\begin{equation}
	\label{eq:Z_tau_t_and_z_t}
	E[||z^{\bar{\pi}}_\e(\tau_t)-z_t||] \le M\Delta(\e)e^{C}, \ \ \  t=0,1,\ldots,N(\e).
\end{equation}
Combining (\ref{eq:bar_z_tau_t_and_zeta_t}), (\ref{eq:zeta_t_and_z_t}) and (\ref{eq:Z_tau_t_and_z_t})
\begin{equation}
	\begin{split}
		\label{eq:z_tau_t_Z_tau_t_proof_inequality}
		E[||\bar{z}(\tau_t) - z^{\bar{\pi}}_\e(\tau_t)||] &\le E[||\bar{z}(\tau_t) - \zeta_t|| + ||\zeta_t - z_t|| + ||z_t - z^{\bar{\pi}}_\e(\tau_t)||]\\
		&\le CM\Delta(\e) + \frac{1}{C}\left(C^2M\Delta(\e) +c\left(\frac{\e}{\Delta(\e)}\right) +  \bar{\nu}_g(K(\e))+CM\Delta(\e)\right)e^{C},
	\end{split}
\end{equation}
for all $t=0,1,\ldots, N(\e)$. Finally, by Lemma \ref{lem:Lipschitz_deterministic_system}, Lemma \ref{lem:hybrid_system_soln_tightness} and (\ref{eq:z_tau_t_Z_tau_t_proof_inequality}),
\begin{equation*}
	\begin{split}
		E||z^{\bar{\pi}}_\e(s) - &\bar{z}(s)|| \le E||\bar{z}(s) - \bar{z}(\tau_t)|| + E[||\bar{z}(\tau_t) - z^{\bar{\pi}}_\e(\tau_t)||] + E||z^{\bar{\pi}}_\e(s) - z^{\bar{\pi}}_\e(\tau_t)|| \\
		& \le 2M\Delta(\e) + CM\Delta(\e) + \frac{1}{C}\left(C^2M\Delta(\e) +c\left(\frac{\e}{\Delta(\e)}\right) +  \bar{\nu}_g(K(\e))+CM\Delta(\e)\right)e^{C}.
	\end{split}
\end{equation*}
The function $\psi_2(\e) \BYDEF M\Delta(\e) + CM\Delta(\e) + \frac{1}{C}\left(C^2M\Delta(\e) +c\left(\frac{\e}{\Delta(\e)}\right) +  \bar{\nu}_g(K(\e))+CM\Delta(\e)\right)e^{C}$ satisfies the conditions of the $\psi_2(\e)$ proposed in the statement of the Theorem, thus, the proof is complete.  $\ \Box$

\section{Proofs of Theorems \ref{thm:expected_occ_set_and_V_h_distance}, \ref{thm:LOMS_theorem_2} and \ref{thm:key_theorem_from_borkar}}
\label{sec:borkar_proofs_1}
\begin{Lemma}
	Let the vector function $h(\cdot,\cdot)$ be as in (\ref{eq:h-AC_function}) and define the constant $c_h$ as
	\begin{equation}\label{def:c_h}
		c_h \BYDEF \max_{(y,u) \in Y \times \hat{U}} ||h(y,u)||.
	\end{equation}
	Then, the following estimates are valid for any pair of integers $T', T'' \ge 1$,
	\begin{align}
		\label{eq:E_zeta_inequality}
		\sup_{\zeta \in V_h^r(T,y_0)} E[||\zeta||] \le c_h \implies \sup_{\zeta \in V_h(T,y_0)}||\zeta|| \le c_h \\
		\label{eq:d_E_H_V_T'_T''_inequality}
		d^E_H\left(V_h^r(T',y_0), V_h^r(T'',y_0)\right) \le \frac{2c_h|T'' - T'|}{\max\{T',T''\}} \\
		\label{eq:d_H_E_V_T'_T''_inequality}
		d_H(V_h(T',y_0),V_h(T'',y_0)) \le \frac{2c_h|T'' - T'|}{\max\{T',T''\}}.
	\end{align}

\end{Lemma}
{\bf Proof.}
The first inequality in (\ref{eq:E_zeta_inequality}) follows from (\ref{def:c_h}), and the second inequality in (\ref{eq:E_zeta_inequality}) is implied by the fact $\|E[\zeta]\| \le E[\|\zeta\|]$. Equation (\ref{eq:d_H_E_V_T'_T''_inequality}) follows from (\ref{eq:d_E_H_V_T'_T''_inequality}), since $d_H(E[V_1],E[V_2])$ $\le$ $d_H^E (V_1, V_2)$. It remains to show (\ref{eq:d_E_H_V_T'_T''_inequality}) holds.
Assume without loss of generality that $T'' \ge T'$. Then, by (\ref{eq:E_zeta_inequality}), for any control plan $\pi \in \Pi$ with the corresponding state-control trajectory $(y^{\pi, y_0}(\cdot),u^{\pi, y_0}(\cdot))$, we have 
\begin{equation*}
	\begin{split}
		E\left[\left|\left| \frac{1}{T'}\sum_{t= 0}^{T'-1}h(y^{\pi, y_0}(t), u^{\pi, y_0}(t)) - \frac{1}{T''}\sum_{t= 0 }^{T'' - 1}h(y^{\pi, y_0}(t), u^{\pi, y_0}(t)) \right|\right|\right] \\
		= E\left[\left|\left| \left(\frac{1}{T'}-\frac{1}{T''}\right)\sum_{t= 0}^{T'-1}h(y^{\pi, y_0}(t), u^{\pi, y_0}(t)) - \frac{1}{T''}\sum_{t= T' }^{T'' - 1}h(y^{\pi, y_0}(t), u^{\pi, y_0}(t)) \right|\right|\right] \\
		\le E\left[\left|\left| \left(\frac{T'' - T'}{T''}\right)\frac{1}{T'}\sum_{t= 0}^{T'-1}h(y^{\pi, y_0}(t), u^{\pi, y_0}(t)) \right|\right|   + \left|\left| \frac{1}{T''}\sum_{t= T'}^{T''-1}h(y^{\pi, y_0}(t), u^{\pi, y_0}(t)) \right|\right|\right]\\
		\le \frac{2c_h (T'' - T')}{T''}
	\end{split}
\end{equation*}
which implies (\ref{eq:d_E_H_V_T'_T''_inequality}).
$\ \Box$
\begin{Lemma}
	The strong $h$-approximation conditions are equivalent to the fulfilment of the inequality
	\begin{equation}\label{eq:strong_h_approx_lemma}
		d^E_H(V_h^r(T,y_0'),V_h^r(T,y_0'')) \le \nu_h(T)
	\end{equation}
	and the weak $h$-approximation conditions are equivalent to the fulfilment of the inequality
	\begin{equation}\label{eq:weak_h_approx_lemma}
		d_H(V_h(T,y_0'),V_h(T,y_0'')) \le \nu_h(T)
	\end{equation}
	for any initial conditions $y_0'$ and $y_0''$.
\end{Lemma}
{\bf Proof.} The proof is obvious. $\ \ \Box$ 
\begin{Lemma}
	\label{lem:lemma_6_1}
	Let a function $\ph(T):\mathbb{N}\rightarrow \mathbb{R}$, where $\mathbb{N}$ denotes the positive integers, be such that, for some monotone decreasing function $\beta(T)$, $\lim\limits_{T\rightarrow \infty} \beta(T) = 0$, the following inequalities hold:
	\begin{equation}
		\label{eq:T_kT_inequality}
		|\ph(T) - \ph(kT)| \le \beta(T), \quad k=1,2,\ldots.
	\end{equation}
	Let also
	\begin{equation}
		\label{eq:T''_T'_inequality}
		|\ph(T') - \ph(T'')| \le \frac{\alpha|T'' - T'|}{\max\{T', T''\}} \quad \forall \ T', T'' > 0, \text{   } \alpha \text{ const.}
	\end{equation} 
	Then, there exists a limit
	\begin{equation}
		\label{eq:A_limit}
		\lim_{T \rightarrow \infty}\ph(T) \BYDEF A
	\end{equation}
	and the estimate
	\begin{equation}
		\label{eq:A_estimate}
		|\ph(T) - A| \le \beta(T) \quad  \forall \ T > 0
	\end{equation}
	is valid.
\end{Lemma}

{\bf Proof.} Let $\delta>0$ be arbitrarily small. Let $T_\delta$ be such that $\beta(T_\delta)\leq \frac{\delta}{4}$ and let $T_\delta^+$  be such that $T_\delta^+\geq \frac{4 \alpha T_\delta}{\delta}$. Then, for 
$T''\geq T'\geq \max\{ T_\delta , T_\delta ^+\}$, we obtain:
\begin{equation*}
	\begin{split}
		\left|\varphi(T'')-\varphi(T')\right|\leq &\left|\varphi(T'')-\varphi\left(\left\lfloor\frac{T''}{T_\delta}\right\rfloor T_\delta\right)\right|+\left|\varphi\left(\left\lfloor\frac{T''}{T_\delta}\right\rfloor T_\delta\right)-\varphi(T_\delta)\right|\\&+\left|\varphi\left(\left\lfloor\frac{T'}{T_\delta}\right\rfloor T_\delta\right)-\varphi(T_\delta)\right|+\left|\varphi(T')-\varphi\left(\left\lfloor\frac{T'}{T_\delta}\right\rfloor T_\delta\right)\right|\\
		&\leq \frac{\alpha{T_\delta}}{T''}+2\beta(T_\delta)+\frac{\alpha{T_\delta}}{T'}\leq 2\beta(T_\delta)+\frac{2\alpha{T_\delta}}{T_\delta^+}\leq \delta.
	\end{split}
\end{equation*} Equation (\ref{eq:A_limit}) follows and (\ref{eq:A_estimate}) is obtained by taking the limit as $k\rightarrow\infty$ in (\ref{eq:T_kT_inequality}). $\ \Box$

\begin{Lemma}
	Let $\Psi(p,T,y_0)$ be the support function of the set $V_h(T,y_0)$:
	\begin{equation}
		\Psi_h(p,T,y_0) \BYDEF \sup_{v \in V_h(T,y_0)} \{p^Tv\}.
	\end{equation}
	If the weak $h$-approximation conditions are satisfied for any vector function $h(\cdot,\cdot)$ as in (\ref{eq:h-AC_function}), then there exists a convex, positively homogeneous and Lipschitz continuous function $\Psi_h(p)$ such that
	\begin{equation}
		\label{eq:lem_6.4}
		|\Psi_h(p,T,y_0) - \Psi_h(p)| \le \nu_h(T)||p||,
	\end{equation}
	where $\nu_h(T)$ is the function introduced in Definition \ref{def:S-h-AC}.
\end{Lemma}
{\bf Proof.}
Note, first, that $\Psi_h(p,T,y_0)$ also allows the representation
\begin{equation*}
	\Psi_h(p,T,y_0) = \frac{1}{T}\sup_{\pi \in \Pi}\left\{ E\left[ \sum_{\tau = 0}^{T - 1}p^Th(y^{\pi, y_0}(\tau),u^{\pi, y_0}(\tau)) \right] \right\}
\end{equation*}
where the supremum is over all control plans $\pi \in \Pi$. From (\ref{eq:E_zeta_inequality}) it follows that 
\begin{equation*}
	|\Psi_h(p,T,y_0)| \le c_h||p||, \quad |\Psi_h(p', T, y_0) - \Psi_h(p'', T, y_0)| \le c_h||p'-p''||
\end{equation*}
and from (\ref{eq:d_H_E_V_T'_T''_inequality}) it follows that
\begin{equation*}
	\begin{split}
		|\Psi_h(p, T', y_0) - \Psi_h(p, T'', y_0)| &=  ||p|| \left|\Psi_h\left(\frac{p}{||p||}, T', y_0\right) - \Psi_h\left(\frac{p}{||p||}, T'', y_0\right) \right|\\
		&\le ||p||\max_{||p'|| \le 1} \left|\Psi_h\left(p', T', y_0\right) - \Psi_h\left(p', T'', y_0\right) \right|\\
		&\le ||p|| d_H(V_h(T',y_0)], V_h(T'',y_0)])\\
		&\le \frac{2c_h||p|||T'' - T'|}{\max\{T',T''\}}
	\end{split}
\end{equation*}
similarly, by (\ref{eq:weak_h_approx_lemma}),
\begin{equation}
	\label{eq:Psi_h_approximation}
	|\Psi_h(p, T, y_0') - \Psi_h(p, T, y_0'')| \le ||p|| \nu_h(T).
\end{equation}
Note that if (\ref{eq:lem_6.4}) is established, then the fact that $\Psi_h(p)$ is convex, positively homogeneous and Lipschitz continuous will follow from the fact that $\Psi_h(p,T,y_0)$ is. By Lemma \ref{lem:lemma_6_1}, to establish (\ref{eq:lem_6.4}) it is sufficient to verify the validity of the following estimates:
\begin{equation}
	\label{eq:Psi_h_kT_and_T}
	|\Psi_h(p,kT,y_0) - \Psi_h(p,T,y_0) | \le c\nu_h(T) \quad k = 1, 2, \ldots.
\end{equation}
For $k=1$, it is obvious. Assume that
\begin{equation}
	\label{eq:Psi_inductive_assumption}
	|\Psi(p, (k-1)T, y_0) - \Psi(p, T, y_0)| \le c\nu_h(T)
\end{equation}
and show the validity of (\ref{eq:Psi_h_kT_and_T}) using induction. Define the collection of random variables $\mathcal{W}_h(T, y_0)$ as follows:
\begin{equation*}
	\mathcal{W}_h(T,y_0) \BYDEF \bigcup_{\pi \in \Pi}  \left\{\left(\frac{1}{T}\sum_{\tau = 0 }^{T-1}h(y^{\pi, y_0}(\tau),u^{\pi, y_0}(\tau)), y^{\pi, y_0}(T)\right)\right\}.
\end{equation*}
Using dynamic programming, one can obtain
\begin{equation}
	\label{eq:Psi_dynamic_prog}
	\Psi_h(p, kT, y_0) = \sup_{(\zeta, \eta) \in \mathcal{W}_h((k-1)T, y_0)} \left\{\frac{k-1}{k}E[p^T\zeta] + \frac{1}{k}E[\Psi_h(p,T,\eta)]  \right\}.
\end{equation}
By (\ref{eq:Psi_h_approximation}),  for any $(\zeta, \eta) \in \mathcal{W}_h((k-1)T, y_0)$, 
\begin{equation*}
	|E[\Psi_h(p,T,\eta)] - \Psi_h(p,T,y_0) | \le  ||p|| \nu_h(T)
\end{equation*}
Hence, using (\ref{eq:Psi_dynamic_prog}), one can obtain
\begin{equation*}
	\begin{split}
		\left|\Psi_h(p,kT,y_0) - \left(\frac{k-1}{k}\Psi_h(p,(k-1)T, y_0) + \frac{1}{k}\Psi_h(p,T,y_0)\right)  \right|\\
		= \left|\Psi_h(p,kT,y_0) - \sup_{(\zeta, \eta) \in \mathcal{W}_h((k-1)T, y_0)}\left\{\frac{k-1}{k}E[p^T\zeta] + \frac{1}{k}\Psi_h(p,T,y_0)\right\}  \right| \\
		\le \frac{1}{k}\sup_{(\zeta, \eta) \in \mathcal{W}_h((k-1)T, y_0)}\{|E[\Psi_h(p,T,\eta)] - \Psi_h(p,T,y_0) |\} \le \left(\frac{1}{k}\right)||p||\nu_h(T)
	\end{split}
\end{equation*}
From (\ref{eq:Psi_inductive_assumption}) it follows, on the other hand, that 
\begin{equation*}
	\left| \frac{k-1}{k}\Psi_h(p,(k-1)T, y_0) - \frac{k-1}{k}\Psi_h(p,T,y_0)) \right| \le \frac{k-1}{k}||p|| \nu_h(T),
\end{equation*}
which implies
\begin{equation*}
	\begin{split}
		|\Psi_h(p,kT,y_0) - \Psi_h(p,T,y_0)| \\
		\le \left|\Psi_h(p,kT,y_0) - \left(\frac{k-1}{k}\Psi_h(p, (k-1)T, y_0) + \frac{1}{k}\Psi_h(p,T,y_0)\right)\right| \\
		+ \left|\frac{k-1}{k}\Psi_h(p, (k-1)T, y_0) - \frac{k-1}{k}\Psi_h(p,T,y_0)\right| \le \left(\frac{1}{k} + \frac{k-1}{k}\right) ||p|| \nu_h(T),
	\end{split}
\end{equation*}
as required.
$\ \Box$
\begin{Lemma}
	Let the system (\ref{eq:associated_system}) satisfy the the weak $h$-approximation conditions for any vector function $h(\cdot,\cdot)$ as in (\ref{eq:h-AC_function}), and let $V_h$ be a convex and compact subset of $\mathbb{R}^j$ defined by 
	\begin{equation}
		\label{eq:V_h_definition_support}
		V_h \BYDEF \{v | p^Tv \le \Psi_h(p) \quad \forall \ p \in \mathbb{R}^j \}.
	\end{equation}
	Then, for any initial condition $y_0$,
	\begin{equation}
		\label{eq:d_H_co_V_h_and_V_h}
		d_H(\text{co}V_h(T,y_0), V_h) \le \nu_h(T).
	\end{equation}
\end{Lemma}
REMARK. The notation (\ref{eq:V_h_definition_support}) anticipates the fact that this set will coincide with the set $V_h$, the existence of which is claimed by Theorem \ref{thm:expected_occ_set_and_V_h_distance}.\\
{\bf Proof.}
Note the fact that the set $V_h$ is convex and compact follows from its definition in (\ref{eq:V_h_definition_support}).
Note also that the support functions of both $V_h(T,y_0)$ and $coV_h(T,y_0)$ are equal to $\Psi_h(p,T,y_0)$:
\begin{equation*}
	\sup_{v \in coV_h(T,y_0)} \left\{p^T v\right\} = \sup_{v \in V_h(T,y_0)} \left\{p^T v\right\} = \Psi_h(p,T,y_0).
\end{equation*}
The support function for $V_h$ is $\Psi_h(p)$ (see Corollary 13.2.1 in \cite{rockafellar1970convex}). Hence (see e.g. Lemma $\Pi$2.9, p. 207 in \cite{gaitsgory1991slowfast}),
\begin{equation}
	\label{eq:d_H_co_E_V_h_and_V_h}
	d_H(\text{co}V_h(T,y_0), V_h) \le \sup_{||p|| \le 1} |\Psi_h(p,T,y_0) - \Psi_h(p)|.
\end{equation}
By (\ref{eq:lem_6.4}), 
\begin{equation*}
	| \Psi_h(p,T,y_0) - \Psi_h(p)| \le ||p|| \nu_h(T).
\end{equation*}
The above and (\ref{eq:d_H_co_E_V_h_and_V_h}) imply (\ref{eq:d_H_co_V_h_and_V_h}).
$\ \Box$
\begin{Proposition}\label{prop:borkar_prop_3-5}
	Let $V_i, i=1,\ldots, k$ be collections of random variables defined on the same probability space such that any element $\zeta_i \in V_i$ is independent from any other element $\zeta_j \in V_j$ for $i \ne j$. Assume also that
	\begin{equation}
		E[||\zeta||^2]\le \bar{c} = const \ \ \forall \zeta \in V_i, \ \ i=1,2,\ldots, k.
	\end{equation}
	Then
	\begin{equation}\label{eq:Vi_EVi_HD}
		d^E_H\left(\frac{1}{k}\sum_{i=1}^{k}V_i, \frac{1}{k}\sum_{i=1}^{k}E[V_i]\right) \le \sqrt{\frac{\bar{c}}{k}}
	\end{equation}
	where $E[V_i]$ stands for the set of mathematical expectations of the elements of the set $V_i$ and $d^E_H(\cdot,\cdot)$ is as defined in (\ref{def:expected_hausdorff_distance}).
\end{Proposition}
{\bf Proof.} Take an arbitrary element $\zeta \in \frac{1}{k}\sum_{i=1}^{k}V_i$. By definition, $\zeta$ may be presented in the form $\zeta = \frac{1}{k}\sum_{i=1}^{k}\zeta_i$, where $\zeta_i \in V_i$. Define
\begin{equation*}
	\bar{\zeta} \BYDEF E[\zeta]= \frac{1}{k}\sum_{i=1}^{k}E[\zeta_i] \in \frac{1}{k}\sum_{i=1}^{k}E[V_i].
\end{equation*}
Due to the independence of $\zeta_i$, $i=1,2,\ldots,k$, and by (\ref{eq:Vi_EVi_HD}),
\begin{equation}\label{exp_rv_rv_distance_1}
	E[||\zeta - \bar{\zeta}||] \le \sqrt{E[||\zeta - \bar{\zeta} ||^2]} = \sqrt{\frac{1}{k^2}\sum_{i=1}^{k}E[||\zeta_i - E[\zeta_i] ||^2]} \le \sqrt{\frac{\bar{c}}{k}}.
\end{equation}
Now take an arbitrary $\bar{\zeta} \in \frac{1}{k}\sum_{i=1}^{k}E[V_i]$. By definition, there exist $\zeta_i \in V_i$, $i=1,2,\ldots,k$ such that $\bar{\zeta} = \frac{1}{k}\sum_{i=1}^{k}E[\zeta_i]$. Define $\zeta \BYDEF \frac{1}{k}\sum_{i=1}^{k}\zeta_i \in \frac{1}{k}\sum_{i=1}^{k}V_i$. By a similar step to (\ref{exp_rv_rv_distance_1}), one can establish that $E[||\zeta - \bar{\zeta}||] \le \sqrt{\frac{\bar{c}}{k}}$, completing the proof of the proposition. $\ \Box$
%
\begin{Lemma}
	For  any $T = 1,2,\ldots$ and $k=1,2,\ldots,$ there exists a collection of random variables $\bar{V}^r_h(kT, y_0)$ such that
	\begin{equation}\label{eq:vg_semi_markov_lemma_a}
		\bar{V}^r_h(kT, y_0) \subset V_h^r( kT, y_0) \ \implies \ E[\bar{V}^r_h(kT, y_0)] \subset V_h(kT, y_0)
	\end{equation}
	and such that: \\
	\textnormal{(i)} the estimate 
	\begin{equation}\label{eq:vg_semi_markov_lemma_b}
		d_H(E[\bar{V}^r_h(kT, y_0)], coV_h(T,y_0)) \le \frac{\bar{c}_h}{k} + \nu_h(T)
	\end{equation}
	is valid if the weak $h$-approximation conditions are satisfied, and\\
	\textnormal{(ii)} the estimate
	\begin{equation}\label{eq:vg_semi_markov_lemma_c}
		d^E_H(\bar{V}^r_h(kT, y_0), coV_h(T,y_0)) \le \frac{c_h}{\sqrt{k}} + \frac{\bar{c}_h}{k} + \nu_h(T)
	\end{equation}
	is valid if the strong $h$-approximation conditions are satisfied, where $c_h$ and $\bar{c}_h$ are constants.
\end{Lemma}
{\bf Proof.} The proof is broken down into three parts.\\
\textbf{Construction of $\bar{V}^r_h(kT, y_0)$.} For $i=1,\ldots, k$, denote by $\Pi^{(i-1)T,iT}$ the family of control plans defined on the time sequence $t=(i-1)T, (i-1)T+1,\ldots,iT-1$ as follows: $\pi^i \in \Pi^{(i-1)T,iT}$ if and only if
$$\pi^i=\big(\pi^i_{(i-1)T}, \pi^i_{(i-1)T+1}(\sigma((i-1)T)),\ldots, \pi^i_{iT-1}(\sigma((i-1)T), \ldots, \sigma(iT-2))\big),$$
where $\pi^i_{(i-1)T}$ is an element of $\hat U$ and $\pi^i_{(i-1)T+t}$ is a Borel measurable function of $\sigma((i-1)T),\ldots,\sigma((i-1)T+t-1)$ for $1 \le t \le T-1$. For $i=1,2,\ldots,k$, define the collection of random variables $V_h^r( (i-1)T, iT, \eta)$:
\begin{equation}
	V_h^r( (i-1)T, iT, \eta) \BYDEF \bigcup_{ \pi^i \in \Pi^{(i-1)T,iT}}\left\{\frac{1}{T}\sum_{t=(i-1)T}^{iT-1}h(y^{\pi^i, \eta}(t),u^{\pi^i, \eta}(t))  \right\} 
\end{equation}
and the set of corresponding mathematical expectations $E[V_h^r( (i-1)T, iT, \eta)]$:
\begin{equation}
	E[V_h^r( (i-1)T, iT, \eta)] \BYDEF \bigcup_{ \pi^i \in \Pi^{(i-1)T,iT}}\left\{E\left[\frac{1}{T}\sum_{t=(i-1)T}^{iT-1}h(y^{\pi^i, \eta}(t),u^{\pi^i, \eta}(t)) \right] \right\},
\end{equation}
where $(y^{\pi^i, \eta}(\cdot), u^{\pi^i, \eta}(\cdot))$ is the state-control trajectory corresponding to a control $\pi^i \in \Pi^{(i-1)T, iT}$ such that $y^{\pi^i, \eta}((i-1)T)= \eta$ is satisfied, where $\eta$ is a random variable independent of the random variables $\sigma(t)$ for $t \ge (i-1)T$. Also, let 
\begin{equation*}
	V_h^r(0,T,y_0) \BYDEF V_h^r(T,y_0),
\end{equation*}
\begin{equation*}
	E[V_h^r(0,T,y_0)] \BYDEF V_h(T,y_0).
\end{equation*}
Denote by $\Pi^{0,kT}$ a family of truncated control plans defined in the following way: a control plan $\pi$ belongs to the set $\Pi^{0,kT}$ if and only if there exists a sequence of control plans $\pi^i \in \Pi^{(i-1)T,iT}$, $i=1,2,\ldots, k$, such that 
\begin{equation*}
	\begin{split}
		\pi = (&\pi^1_0, \pi^1_1(\sigma(0)),\ldots,\pi^1_{T-1}(\sigma(0),\ldots,\sigma(T-2)), \ldots, \\  
		&\pi^i_{(i-1)T}, \pi^i_{(i-1)T+1}(\sigma((i-1)T)),\ldots,\pi^i_{iT-1}(\sigma((i-1)T),\ldots,\sigma(iT-2)),\ldots,\\ &\pi^k_{(k-1)T}, \pi^k_{(k-1)T+1}(\sigma((k-1)T)),\ldots,\pi^k_{kT-1}(\sigma((k-1)T),\ldots,\sigma(kT-2))).
	\end{split}
\end{equation*}

Define $\bar{V}^r_h(kT, y_0)$ as the collection of random variables
\begin{equation}
	\bar{V}^r_h(kT, y_0) \BYDEF \bigcup_{ \pi \in \Pi^{0,kT} }\left\{\frac{1}{kT}\sum_{t=0}^{kT-1}h( y^{\pi, y_0}(t),u^{\pi, y_0}(t)) \right\}
\end{equation}
and $E[\bar{V}^r_h(kT, y_0)]$ as the set of corresponding mathematical expectations
\begin{equation}
	E[\bar{V}^r_h(kT, y_0)] \BYDEF \bigcup_{\pi \in \Pi^{0,kT}}\left\{E\left[\frac{1}{kT}\sum_{t=0}^{kT-1}h( y^{\pi, y_0}(t),u^{\pi, y_0}(t))\right]  \right\}
\end{equation}
where, in both cases, the union is over the control plans $\pi \in \Pi^{0,kT}$ with the corresponding state-control trajectories $(y^{\pi, y_0}(\cdot), u^{\pi, y_0}(\cdot))$. Note that, by definition,
\begin{equation}
	\bar{V}^r_h(T,y_0) = V_h^r(T,y_0), \  E[\bar{V}^r_h(T,y_0)] = V_h(T,y_0),
\end{equation}
and that inclusions (\ref{eq:vg_semi_markov_lemma_a}) are valid for $k=2,3,\ldots$. It may be readily verified that
\begin{equation}\label{eq:vh'_dyn_prog}
	\bar{V}^r_h(kT,y_0) = \bigcup_{(\zeta, \eta) \in \mathcal{W}_h'((k-1)T,y_0)}\left\{\frac{k-1}{k}\zeta + \frac{1}{k}V_h^r((k-1)T,kT,\eta) \right\}
\end{equation}
and
\begin{equation}\label{eq:vh'_expectation_dyn_prog}
	E[\bar{V}^r_h(kT,y_0)] = \bigcup_{(\zeta, \eta) \in \mathcal{W}_h'( (k-1)T,y_0)}\left\{\frac{k-1}{k}E[\zeta] + \frac{1}{k}E[V_h^r((k-1)T,kT,\eta)] \right\}
\end{equation}
where
\begin{equation}
	\mathcal{W}_h'( (k-1)T,y_0) = \bigcup_{\pi \in \Pi^{0,(k-1)T}} \left\{\left(\frac{1}{(k-1)T}\sum_{t=0}^{(k-1)T-1}h(y^{\pi, y_0}(t),u^{\pi, y_0}(t)),	y^{\pi, y_0}((k-1)T)  \right)  \right\},
\end{equation}
the union being over the controls plans from $\Pi_0^{(k-1)T}$.
\\
\textbf{Proof of (i).} Using induction, let us show that
\begin{equation}\label{eq:vh_lemma_inductive_proof_a}
	d_H(E[\bar{V}^r_h(kT, y_0)],\frac{1}{k}\sum_{i=1}^{k}E[V_h^r((i-1)T,iT,y_0)]) \le \nu_h(T), \ \ k=1,2,\ldots
\end{equation}
For $k=1$ it is immediate since, by definition, $E[\bar{V}^r_h(T,y_0)]=E[V_h^r(0,T,y_0)]$. Assume that the estimate 
\begin{equation}\label{eq:vh_lemma_inductive_assumption_a}
	d_H(E[\bar{V}^r_h((k-1)T, y_0)],\frac{1}{k-1}\sum_{i=1}^{k}E[V_h^r((i-1)T,iT,y_0)]) \le \nu_h(T), \ \ k=1,2,\ldots
\end{equation}
is valid. From the weak $h$-approximation conditions (see (\ref{eq:weak_h_approx_lemma})) it follows that, for any $\eta$ such that $(\zeta, \eta) \in \mathcal{W}_h'((k-1)T, y_0)$,
\begin{equation}
	d_H(E[V_h^r((i-1)T,iT,\eta)]), E[V_h^r((k-1)T,iT,y_0)]) \le \nu_h(T).
\end{equation}
This and (\ref{eq:vh'_expectation_dyn_prog}) lead to the estimate
\begin{equation}
	\begin{split}
		&d_H(E[\bar{V}^r_h(kT, y_0)], \frac{k-1}{k}E[\bar{V}^r_h((k-1)T, y_0)]+\frac{1}{k}E[V_h^r((k-1)T,kT,y_0)]) \\
		&=d_H\left(E[\bar{V}^r_h(kT, y_0)], \bigcup_{(\zeta,\eta)\in\mathcal{W}_h'((k-1)T,y_0)}\left\{\frac{k-1}{k}E[\zeta]+\frac{1}{k}E[V_h^r((k-1)T,kT,y_0)]\right\}\right) \\
		&\le d_H\left(E[\bar{V}^r_h(kT, y_0)], \bigcup_{(\zeta,\eta)\in\mathcal{W}_h'((k-1)T,y_0)}\left\{\frac{k-1}{k}E[\zeta]+\frac{1}{k}E[V_h^r((k-1)T,kT,\eta)]\right\}\right) + \left(\frac{1}{k}\right)\nu_h(T) \\
		&= \left(\frac{1}{k}\right)\nu_h(T).
	\end{split}
\end{equation}
Using the estimate above and (\ref{eq:vh_lemma_inductive_assumption_a}), one can further obtain that
\begin{equation}
	\begin{split}
		&d_H\left(E[\bar{V}^r_h(kT,y_0)],\frac{1}{k}\sum_{i=1}^{k}E[V_h^r((i-1)T,iT,y_0)]\right) \\
		& \le d_H\left(\frac{k-1}{k}E[\bar{V}^r_h((k-1)T, y_0)]+\frac{1}{k}E[V_h^r((k-1)T,kT,y_0)], \frac{1}{k}\sum_{i=1}^{k}E[V_h^r((i-1)T,iT,y_0)]\right) \\ & \ \ \ + \left(\frac{1}{k}\right)\nu_h(T)\\
		& = d_H\left(\frac{k-1}{k}E[\bar{V}^r_h((k-1)T,y_0)], \frac{k-1}{k}\frac{1}{k-1}\sum_{i=1}^{k-1}E[V_h^r((i-1)T,iT,y_0)]\right)+ \left(\frac{1}{k}\right)\nu_h(T) \\
		& \le \left(\frac{k-1}{k}\right)\nu_h(T) + \left(\frac{1}{k}\right)\nu_h(T)\\
		&= \nu_h(T).
	\end{split}
\end{equation}
Thus, (\ref{eq:vh_lemma_inductive_proof_a}) has been established. Since 
\begin{equation}
	E[V_h^r((i-1)T,iT,y_0)] = V_h(T,y_0) \ \ \forall \ i=1,\ldots,k,
\end{equation}
equation (\ref{eq:vh_lemma_inductive_proof_a}) is equivalent to
\begin{equation}\label{eq:vh_lemma_weak_triangle_a}
	d_H(E[\bar{V}^r_h(kT, y_0)],\frac{1}{k}\sum_{i=1}^{k}V_h(T,y_0)) \le \nu_h(T), \ \ k=1,2,\ldots
\end{equation}
By Shapley-Folkman's theorem (see, e.g., \cite[p. 220]{gaitsgory1991slowfast}), 
\begin{equation}\label{eq:vh_lemma_weak_triangle_b}
	d_H\left(\frac{1}{k}\sum_{i=1}^{k}V_h(T,y_0), coV_h(T,y_0)\right) \le \frac{2(j+1)c_h}{k}
\end{equation}
where $c_h$ is as defined in (\ref{def:c_h}) and $j$ is the number of the components of the vector function $h(\cdot,\cdot)$ (see (\ref{eq:h-AC_function})). Equations (\ref{eq:vh_lemma_weak_triangle_a}) and (\ref{eq:vh_lemma_weak_triangle_b}) imply (\ref{eq:vg_semi_markov_lemma_b}) with $\bar{c}_h \BYDEF 2(j+1)c_h$.
\\\textbf{Proof of (ii).} 
Using induction, let us show that
\begin{equation}\label{eq:vh_lemma_inductive_proof_b}
	d^E_H\left(\bar{V}^r_h(kT, y_0), \frac{1}{k}\sum_{i=1}^{k}V_h^r((i-1)T, iT, y_0)\right) \le \nu_h(T), \ k=1,2,\ldots.
\end{equation}
For $k=1$, the above expression is obvious. Assume that 
\begin{equation}\label{eq:vh_lemma_inductive_assumption_b}
	d^E_H\left(\bar{V}^r_h((k-1)T, y_0), \frac{1}{k-1}\sum_{i=1}^{k-1}V_h^r((i-1)T, iT, y_0)\right) \le \nu_h(T).
\end{equation}
From the strong $h$-approximation conditions (see (\ref{eq:strong_h_approx_lemma})) it follows that, for any $\eta$ such that $(\zeta, \eta) \in \mathcal{W}_h'((k-1)T, y_0)$,
\begin{equation}
	d^E_H(V_h^r((k-1)T, kT, \eta), V_h^r((k-1)T, kT, y_0)) \le \nu_h(T).
\end{equation}
Hence, by (\ref{eq:vh'_dyn_prog}),
\begin{align*}
	\begin{split}
		&d^E_H\left(\bar{V}^r_h(kT, y_0), \frac{k-1}{k}\bar{V}^r_h((k-1)T, y_0) + \frac{1}{k}V_h^r((k-1)T, kT, y_0)\right)\\
		&= d^E_H\left(\bar{V}^r_h(kT, y_0), \bigcup_{(\zeta, \eta) \in \mathcal{W}_h'((k-1)T, y_0)}\left\{\frac{k-1}{k}\zeta + \frac{1}{k}V_h^r((k-1)T, kT, y_0)\right\}\right) \\
		&\le d^E_H\left(\bar{V}^r_h(kT, y_0), \bigcup_{(\zeta, \eta) \in \mathcal{W}_h'((k-1)T, y_0)}\left\{\frac{k-1}{k}\zeta + \frac{1}{k}V_h^r((k-1)T, kT, \eta)\right\}\right) + \left(\frac{1}{k}\right)\nu_h(T) \\
		& = \left(\frac{1}{k}\right)\nu_h(T).
	\end{split}
\end{align*}
Using this estimate and (\ref{eq:vh_lemma_inductive_assumption_b}), one obtains
\begin{equation*}
	\begin{split}
		&d^E_H\left(\bar{V}^r_h(kT, y_0), \frac{1}{k}\sum_{i=1}^{k}V_h^r((i-1)T, iT, y_0)\right) \\
		&\le d^E_H\left(\frac{k-1}{k}\bar{V}^r_h((k-1)T, y_0) + \frac{1}{k}V_h^r((k-1)T, kT, y_0), \frac{1}{k}\sum_{i=1}^{k}V_h^r((i-1)T, iT, y_0)\right) + \left(\frac{1}{k}\right)\nu_h(T) \\
		&= d^E_H\left(\frac{k-1}{k}\bar{V}^r_h((k-1)T, y_0), \frac{k-1}{k}\frac{1}{k-1}\sum_{i=1}^{k-1}V_h^r((i-1)T, iT, y_0)\right) + \left(\frac{1}{k}\right)\nu_h(T)\\
		&\le  \left(\frac{k-1}{k}\right)\nu_h(T)+ \left(\frac{1}{k}\right)\nu_h(T)\\
		&= \nu_h(T).
	\end{split}
\end{equation*}
This proves the validity of the estimate (\ref{eq:vh_lemma_inductive_proof_b}). Note that for $i,j=1,2,\ldots, k$, the elements of $V_h^r((i-1)T, iT, y_0)$ and $V_h^r((j-1)T, jT, y_0)$ are independent for $i \neq j$, and
\begin{equation*}
	E[||\zeta||^2] \le c_h^2 \ \ \forall \zeta \in V_h^r((i-1)T, iT, y_0), \ \ i=1,2,\ldots,k,
\end{equation*}
where $c_h$ is as defined in (\ref{def:c_h}). Hence, one can use Proposition \ref{prop:borkar_prop_3-5} and the fact $E[V_h^r((i-1)T, iT, y_0)] = V_h(T,y_0)$ to obtain
\begin{equation*}
	\begin{split}
		&d^E_H\left( \frac{1}{k}\sum_{i=1}^{k}V_h^r((i-1)T, iT, y_0), \frac{1}{k}\sum_{i=1}^{k}E[V_h^r((i-1)T, iT, y_0)] \right)\\
		&=d^E_H\left( \frac{1}{k}\sum_{i=1}^{k}V_h^r((i-1)T, iT, y_0), \frac{1}{k}\sum_{i=1}^{k}V_h(T,y_0) \right) \le \frac{c_h}{\sqrt{k}}.
	\end{split}
\end{equation*}
The above estimate along with (\ref{eq:vh_lemma_weak_triangle_b}) and (\ref{eq:vh_lemma_inductive_proof_b}) imply (\ref{eq:vg_semi_markov_lemma_c}). 
$\ \Box$
\begin{Corollary}
	For any initial condition $y_0$,
	\begin{equation}\label{eq:vh'_and_vh_corollary_weak}
		d_H(E[\bar{V}^r_h(kT,y_0)], V_h) \le \frac{\bar{c}_h}{k}+2\nu_h(T)
	\end{equation}
	if the weak $h$-approximation conditions are satisfied, and
	\begin{equation}\label{eq:vh'_and_vh_corollary_strong}
		d^E_H(\bar{V}^r_h(kT,y_0), V_h) \le \frac{c_h}{\sqrt{k}}+\frac{\bar{c}_h}{k}+2\nu_h(T)
	\end{equation}
	if the strong $h$-approximation conditions are satisfied.
\end{Corollary}
{\bf Proof.} The estimates follow from (\ref{eq:d_H_co_V_h_and_V_h}), (\ref{eq:vg_semi_markov_lemma_b}) and (\ref{eq:d_H_co_V_h_and_V_h}), (\ref{eq:vg_semi_markov_lemma_c}), respectively.

\begin{Lemma}
	Assume the system (\ref{eq:associated_system}) satisfies the strong $h$-approximation conditions for any vector function $h(\cdot,\cdot)$ as in (\ref{eq:h-AC_function}). 
	For any initial condition $y_0$,
	\begin{equation}\label{eq:V_h_and_E_V_h^r(T,y_0)_distance}
		\sup_{\zeta \in V_h}d(\zeta, V_h(T,y_0)) \le \nu_h^1(T), \quad \lim_{T \rightarrow \infty}\nu_h^1(T) = 0
	\end{equation}
	and
	\begin{equation}\label{eq:V_h_and_E_V_h^r(T,y_0)_Edistance}
		\sup_{\zeta \in V_h}d^E(\zeta, V_h^r(T,y_0)) \le \nu_h^2(T), \quad \lim_{T \rightarrow \infty}\nu_h^2(T) = 0
	\end{equation}	
\end{Lemma}
{\bf Proof.}
Using (\ref{eq:d_H_E_V_T'_T''_inequality}) with $T'' = T$ and $T' = \lfloor T^{1/2} \rfloor \lfloor T^{1/2} \rfloor$, one can obtain 
\begin{equation}\label{eq:FLSQTC}
	\begin{split}
		d_H(V_h(\lfloor T^{1/2} \rfloor \lfloor T^{1/2} \rfloor, y_0), V_h(T, y_0)) &\le \frac{2c_h|T - \lfloor T^{1/2} \rfloor \lfloor T^{1/2} \rfloor|}{\max\{\lfloor T^{1/2} \rfloor \lfloor T^{1/2} \rfloor, T\}} \\
		&\le \frac{2c_h|T - (T^{1/2}-1)(T^{1/2}-1) |}{T}\\
		&\le \frac{2c_h|2T^{1/2}-1 |}{T}\\
		&\le \frac{4c_h}{T^{1/2}}
	\end{split}
\end{equation}
Hence,
\begin{equation}\label{eq:Vh_and_E_Vh(T,y_0)_prelim_distance}
	\sup_{\zeta \in V_h} d(\zeta, V_h(T, y_0)) \le \sup_{\zeta \in V_h} d(\zeta, V_h(\lfloor T^{1/2}\rfloor \lfloor T^{1/2} \rfloor, y_0)) + \frac{4c_h}{T^{1/2}}
\end{equation}
From (\ref{eq:vg_semi_markov_lemma_a}) and (\ref{eq:vh'_and_vh_corollary_weak}) (with the replacement of $T$ by $\lfloor T^{1/2} \rfloor$ and the replacement of $k$ by $\lfloor T^{1/2} \rfloor$) it follows, on the other hand, that
\begin{equation}
	\begin{split}
		\sup_{\zeta \in V_h} d(\zeta,V_h(\lfloor T^{1/2} \rfloor \lfloor T^{1/2} \rfloor, y_0)]) \le \sup_{\zeta \in V_h} d(\zeta,E[\bar{V}^r_h(\lfloor T^{1/2} \rfloor \lfloor T^{1/2} \rfloor, y_0)])\\
		\le d_H(E[\bar{V}^r_h(\lfloor T^{1/2}\rfloor \lfloor T^{1/2} \rfloor, y_0)], V_h) \le \frac{\bar{c}_h}{\lfloor T^{1/2}\rfloor} + 2\nu_h(\lfloor T^{1/2} \rfloor).
	\end{split}
\end{equation}
This and (\ref{eq:Vh_and_E_Vh(T,y_0)_prelim_distance}) imply (\ref{eq:V_h_and_E_V_h^r(T,y_0)_distance}) with $\nu_h^1(T) \BYDEF \frac{4c_h}{T^{1/2}} + \frac{\bar{c}_h}{\lfloor T^{1/2}\rfloor} + 2\nu_h(\lfloor T^{1/2} \rfloor)$. To establish (\ref{eq:V_h_and_E_V_h^r(T,y_0)_Edistance}), one can use (\ref{eq:d_E_H_V_T'_T''_inequality}) and obtain, similarly to (\ref{eq:FLSQTC}), that
\begin{align}
	&d^E_H(V_h^r(\lfloor T^{1/2} \rfloor \lfloor T^{1/2} \rfloor, y_0), V_h^r(T,y_0)) \le \frac{4c_h}{T^{1/2}}, \\
	\label{eq:Vh_and_E_Vh(T,y_0)_prelim_Edistance} \implies& \sup_{\zeta \in V_h} d^E(\zeta, V_h^r(T,y_0)) \le \sup_{\zeta \in V_h}d^E(\zeta, V_h^r(\lfloor T^{1/2}\rfloor \lfloor T^{1/2} \rfloor,y_0))+\frac{4c_h}{T^{1/2}}.
\end{align}
By (\ref{eq:vg_semi_markov_lemma_a}) and (\ref{eq:vh'_and_vh_corollary_strong}) (with the replacement of $T$ by $\lfloor T^{1/2} \rfloor$ and the replacement of $k$ by $\lfloor T^{1/2} \rfloor$),
\begin{equation}
	\begin{split}
		\sup_{\zeta \in V_h}d^E(\zeta, V_h^r(\lfloor T^{1/2}\rfloor \lfloor T^{1/2} \rfloor,y_0))  \le \sup_{\zeta \in V_h}d^E(\zeta, \bar{V}^r_h(\lfloor T^{1/2}\rfloor \lfloor T^{1/2} \rfloor,y_0))\\
		\le d^E_H(\bar{V}^r_h(\lfloor T^{1/2}\rfloor \lfloor T^{1/2} \rfloor, y_0), V_h) \le \frac{c_h}{\sqrt{\lfloor T^{1/2} \rfloor}} + \frac{\bar{c}_h}{\lfloor T^{1/2}\rfloor} + 2\nu_h(\lfloor T^{1/2} \rfloor).
	\end{split}
\end{equation}
This and (\ref{eq:Vh_and_E_Vh(T,y_0)_prelim_Edistance}) imply (\ref{eq:V_h_and_E_V_h^r(T,y_0)_Edistance}) with $\nu_h^2(T) = \frac{4c_h}{T^{1/2}} + \frac{c_h}{\sqrt{\lfloor T^{1/2} \rfloor}} + \frac{\bar{c}_h}{\lfloor T^{1/2}\rfloor} + 2\nu_h(\lfloor T^{1/2} \rfloor)$. $\ \Box$

\noindent
{\bf Proof of Theorem \ref{thm:expected_occ_set_and_V_h_distance}.} Let us first show that the validity of (\ref{eq:expected_occ_set_and_V_h_distance}) implies the weak $h$-approximation conditions are satisfied. Note that if (\ref{eq:expected_occ_set_and_V_h_distance}) holds, then
\begin{equation*}
	d_H(V_h(T,y'_0), V_h(T,y''_0)) \le d_H(V_h(T,y'_0), V_h) + d_H(V_h, V_h(T,y''_0)) \le 2\nu'_h(T)
\end{equation*}
which, by (\ref{eq:weak_h_approx_lemma}), leads to the fulfilment of the weak $h$-approximation. Let us now show that if the weak $h$-approximation conditions are satisfied, then (\ref{eq:expected_occ_set_and_V_h_distance}) is valid, with $V_h$ being as defined in (\ref{eq:V_h_definition_support}). By (\ref{eq:d_H_co_V_h_and_V_h}),
\begin{equation*}
	\sup_{\zeta \in V_h(T,y_0)}d(\zeta, V_h) \le \sup_{\zeta \in coV_h(T,y_0)}d(\zeta, V_h) \le d_H(coV_h(T,y_0), V_h) \le \nu_h(T).
\end{equation*}
Combining this estimate with (\ref{eq:V_h_and_E_V_h^r(T,y_0)_distance}), one obtains (\ref{eq:expected_occ_set_and_V_h_distance}) with $\nu'_h(T) = \max\{\nu_h(T), \nu^1_h(T)\}$. $\ \Box$

\noindent
{\bf Proof of Theorem \ref{thm:LOMS_theorem_2}.} The proof is broken down into five parts:\\
\textbf{Part I}: The set $\G$ defined in (\ref{eq:LOMS_theorem_2_2}) is compact in the weak$^*$ topology.\\
\textbf{Part II}: 
The set $\G$ is convex.\\
\textbf{Part III}: If the system (\ref{eq:associated_system}) satisfies the weak $h$-approximation conditions (see Definition \ref{def:W-h-AC}) for any vector function $h(\cdot,\cdot)$ as in (\ref{eq:h-AC_function}), then there exists a function $\nu(T)$ tending to zero as $T$ tends to infinity such that 
\begin{equation}\label{eq:LOMS_thm_P3}
	\sup_{\g \in \G_T(y_0)} \rho(\g, \G) \le \nu(T), \ \ \ \ \forall \ y_0 \in Y.
\end{equation}
\textbf{Part IV}: There exists a function $\nu(T)$ tending to zero as $T$ tends to infinity such that 
\begin{equation}\label{eq:LOMS_thm_P4}
	\sup_{\g \in \G} \rho(\g, \G_T(y_0)) \le \nu(T), \ \ \ \ \forall \ y_0 \in Y.
\end{equation}
\textbf{Part V}: If there exists a convex and compact set $\G \in \mathcal{P}(Y \times \hat{U})$ which satisfies (\ref{eq:LOMS_theorem_2_1}) for any initial condition $y_0$, then the weak $h$-approximation is satisfied for any vector function $h(\cdot, \cdot)$ as in (\ref{eq:h-AC_function}) and any initial conditions $y'_0$, $y''_0$. Also, for any $h(\cdot, \cdot)$ as in (\ref{eq:h-AC_function}), the estimate (\ref{eq:expected_occ_set_and_V_h_distance}) is valid with $V_h$ as in (\ref{eq:V_h_G_representation}).\\
\textbf{Proof of Part I}: Since $\mathcal{P}(Y \times \hat{U})$ is weak$^*$ compact (by the Banach-Alaoglu theorem; see, e.g., Theorem 3.5.16 in \cite{Ash}), it is enough to show that $\G$ is closed in the weak$^*$ topology. Let $\g^t \in \G$, $t=1,2,\ldots$ and 
\begin{equation}\label{eq:LOMS_thm_P1_1}
	\lim_{t\rightarrow \infty}\rho(\g^t,\g) = 0.
\end{equation}
By (\ref{eq:LOMS_theorem_2_2}), for any $h(\cdot,\cdot)$ as in (\ref{eq:h-AC_function}),
\begin{equation}
	v_t=\int_{Y \times \hat{U}}^{} h(y,u)\g^t(dy,du) \in V_h.
\end{equation}
From (\ref{eq:LOMS_thm_P1_1}), it follows that 
\begin{equation}
	\lim_{t\rightarrow \infty}v_t = \int_{Y \times \hat{U}}^{} h(y,u)\g(dy,du) \BYDEF v
\end{equation}
Since $V_h$ is closed, $v \in V_h$. This holds for any $h(\cdot, \cdot)$ as in (\ref{eq:h-AC_function}). Hence, by the definition of $\G$ (see (\ref{eq:LOMS_theorem_2_2})), it follows that $\g \in \G$.\\
\textbf{Proof of Part II}: Let $\g^1 \in \G$ and $\g^2 \in \G$. We have
\begin{equation}
	v_i \BYDEF \int_{Y \times \hat{U}}^{} h(y,u) \g^i (dy,du) \in V_h \ \ i=1,2,
\end{equation}
for any function $h(\cdot, \cdot)$ as in (\ref{eq:h-AC_function}). Since $V_h$ is convex, for any $\lambda \in [0,1]$,
\begin{equation}
	\int_{Y \times \hat{U}}^{} h(y,u) (\lambda \g^1 + (1-\lambda) \g^2)(dy,du) = \lambda v_1 + (1-\lambda) v_2 \in V_h.
\end{equation}
The above holds for any $h(\cdot, \cdot)$ as in (\ref{eq:h-AC_function}), implying $\lambda \g^1 + (1-\lambda)\g^2 \in \G$.\\
\textbf{Proof of Part III}: Assume (\ref{eq:LOMS_thm_P3}) does not hold. Then there exists a positive number $\a$ and sequences $\{T^i\}, \{y_0^i\}, \{\g^i \}$, $i=1,2,\ldots$ such that 
\begin{equation}\label{eq:LOMS_thm_P3_sequences}
	\lim_{i \rightarrow \infty} T^i = \infty, \ \ y_0^i \in Y, \ \ \g^i \in \G_{T^i}(y_0^i),
\end{equation}
and such that 
\begin{equation}\label{eq:LOMS_thm_P3_2}
	\rho(\g^i, \G) \ge \a > 0, \ \ \ \forall \ i =1,2,\ldots
\end{equation}
Without loss of generality, one can assume that there exists $\g^* \in \mathcal{P}(Y \times \hat{U})$ such that
\begin{equation}\label{eq:LOMS_thm_P3_3}
	\lim_{i \rightarrow \infty}\rho(\g^i, \g^*) = 0.
\end{equation}	
Since $\rho(\cdot,\cdot)$ is a continuous function of its arguments, (\ref{eq:LOMS_thm_P3_2}) and (\ref{eq:LOMS_thm_P3_3}) imply
\begin{equation}
	\rho(\g^*, \G) \ge \a \implies \g^* \notin \G.
\end{equation}
By (\ref{eq:LOMS_theorem_2_2}), the above means that there exists a function $h(\cdot, \cdot)$ as in (\ref{eq:h-AC_function}) such that 
\begin{equation}\label{eq:LOMS_thm_P3_4}
	\int_{Y \times \hat{U}}^{} h(y,u) \g^* (dy,du) \notin V_h.
\end{equation}
However, 
\begin{equation}\label{eq:LOMS_thm_P3_5}
	v_i \BYDEF \int_{Y \times \hat{U}}^{} h(y,u) \g^i (dy,du) \in V_h(T^i ,y^i_0)
\end{equation}
and, by (\ref{eq:LOMS_thm_P3_3}), 
\begin{equation}\label{eq:LOMS_thm_P3_6}
	\lim_{i \rightarrow \infty} v_i = v^* \BYDEF \int_{Y \times \hat{U}}^{} h(y,u) \g^* (dy,du). 
\end{equation}
From (\ref{eq:LOMS_thm_P3_5}), (\ref{eq:LOMS_thm_P3_6}) and by (\ref{eq:expected_occ_set_and_V_h_distance}) it follows that
\begin{equation*}
	v^* \in V_h,
\end{equation*}
which contradicts (\ref{eq:LOMS_thm_P3_4}). Thus, the number $\a$ and sequences (\ref{eq:LOMS_thm_P3_sequences}) do not exist and (\ref{eq:LOMS_thm_P3}) is established.\\
\textbf{Proof of Part IV}: Assume (\ref{eq:LOMS_thm_P4}) does not hold. Then, there exists a positive number $\a$ and sequences $\{T^i\}, \{y_0^i\}, \{\g^i \}$, $i=1,2,\ldots$ such that 
\begin{equation}\label{eq:LOMS_thm_P4_sequences}
	\lim_{i \rightarrow \infty} T^i = \infty, \ \ y_0^i \in Y, \ \ \g^i \in \G.
\end{equation}
and such that
\begin{equation}
	\rho(\g^i, \G_{T^i}(y^i_0)) \ge \a > 0 \ \ \ \ \forall \ i = 1,2,\ldots
\end{equation}
Since $\G$ is compact, one can choose sequences (\ref{eq:LOMS_thm_P4_sequences}) in such a way that for some $\g^* \in \G$, 
\begin{equation}
	\lim_{i \rightarrow \infty}\rho(\g^i,\g^*) =0.
\end{equation}
This implies that for $i \ge i_0$ ($i_0$ is some positive integer)
\begin{equation*}
	\rho(\g^*, \G_{T^i}(y^i_0)) \ge \frac{\a}{2} > 0
\end{equation*}
or 
\begin{equation*}
	\rho(\g^*, \g) \ge \frac{\a}{2}>0, \ \ \ \ \forall \g \  \in \G_{T^i}( y^i_0).
\end{equation*}
By (\ref{e-rho}), the latter is equivalent to
\begin{equation}
	\sum_{j=1}^{\infty} {1\over 2^j}\left|\int_{Y \times \hat{U}} q_j(y,u)\g^*(dy,du)-\int_{Y \times \hat{U}} q_j(y,u)\g(dy,du)\right| \ge \frac{\a}{2}>0, \ \ \ \ \forall \ \g \in \G_{T^i}(y^i_0).
\end{equation}
Hence, for some integer $N > 0$,
\begin{equation}
	\sum_{j=1}^{N} {1\over 2^j}\left|\int_{Y \times \hat{U}} q_j(y,u)\g^*(dy,du)-\int_{Y \times \hat{U}} q_j(y,u)\g(dy,du)\right| \ge \frac{\a}{4}>0, \ \ \ \ \forall \ \g \in \G_{T^i}(y^i_0).
\end{equation}
implying
\begin{equation}\label{eq:LOMS_thm_P4_1}
	\sum_{j=1}^{N}\sqrt{\left|\int_{Y \times \hat{U}} q_j(y,u)\g^*(dy,du)-\int_{Y \times \hat{U}} q_j(y,u)\g(dy,du)\right|^2} \ge \frac{c_N\a}{4}>0, \ \ \ \ \forall \ \g \in \G_{T^i}(y^i_0),
\end{equation}
where $c_N$ is an appropriately chosen constant. Let $h: Y \times \hat{U} \rightarrow \mathbb{R}^N$ be the vector function with the coordinates $q_t$. That is,
\begin{equation}
	h(y,u) = \{q_t(y,u) \}, \ \ \ \ t=1,\ldots,N.
\end{equation}
Then, by (\ref{eq:LOMS_thm_P4_1}) and (\ref{def:distance}),
\begin{equation}
	d\left(\int_{Y \times \hat{U}}^{} h(y,u)\g^* (dy,du), \int_{Y \times \hat{U}}^{} h(y,u)\g(dy,du)\right) \ge \frac{\a}{4}>0 \ \ \ \ \forall \ \g \in \G_{T^i}(y^i_0),
\end{equation}
implying
\begin{equation}
	\inf_{\g \in \G_{T^i}(y^i_0)}d\left(\int_{Y \times \hat{U}}^{} h(y,u)\g^* (dy,du), \int_{Y \times \hat{U}}^{} h(y,u)\g(dy,du)\right) \ge \frac{\a}{4}>0.
\end{equation}
By definition of $V_h^r(T,y_0)$, the above can be re-written
\begin{equation}
	\begin{split}
		\inf_{v \in V_h(T^i, y^i_0)}d\left(\int_{Y \times \hat{U}}^{} h(y,u)\g^* (dy,du), v\right) \ge \frac{\a}{4}>0 \\
		\implies d\left(\int_{Y \times \hat{U}}^{} h(y,u)\g^* (dy,du), V_h(T^i, y^i_0)\right) \ge \frac{\a}{4}>0.
	\end{split}
\end{equation}
Hence, by (\ref{eq:expected_occ_set_and_V_h_distance}) and the above,
\begin{equation}
	\implies d\left(\int_{Y \times \hat{U}}^{} h(y,u)\g^* (dy,du), V_h\right) \ge \frac{\a}{4} - \nu_h'(T)\ge \frac{\a}{8}, \ \ \ \forall \ i > i_1,
\end{equation}
where $i_1$ is some positive integer. Hence,
\begin{equation}
	\int_{Y \times \hat{U}}^{} h(y,u)\g^* (dy,du) \notin V_h,
\end{equation} 
which, by (\ref{eq:LOMS_theorem_2_2}), contradicts the fact $\g^* \in \G$. Thus, the number $\a$ and sequences (\ref{eq:LOMS_thm_P4_sequences}) do not exist and (\ref{eq:LOMS_thm_P4}) is established.\\
\textbf{Proof of Part V}: By definition of $V_h(T, y_0)$ (see (\ref{eq:V_h_T_definition})), it is straightforward to
verify that the validity of (\ref{eq:LOMS_theorem_2_1}) implies the validity of (\ref{eq:expected_occ_set_and_V_h_distance}) with $V_h$ as in (\ref{eq:V_h_G_representation}) for
any $h(\cdot, \cdot)$ as in (\ref{eq:h-AC_function}). The fact that weak $h$-approximation conditions are satisfied for any such $h(\cdot, \cdot)$ follows now
from Theorem \ref{thm:expected_occ_set_and_V_h_distance}. $\ \Box$

\noindent
{\bf Proof of Theorem \ref{thm:key_theorem_from_borkar}.}
Let us first prove that
\begin{equation}\label{eq:main_borkar_result_1}
	\sup_{\g \in \G} \rho^E(\g,B_T(y_0)) \le \bar{\nu}(T), \ \ \ \ \lim_{T \rightarrow \infty}\bar{\nu}(T) = 0,
\end{equation}
Assume (\ref{eq:main_borkar_result_1}) does not hold. Then, there exists a number $\delta > 0$ and sequences $\mu_t \in \G$ and $T_t$, $t=1,2,\ldots$ (where $T_t \rightarrow \infty$ as $t \rightarrow \infty$) such that, for any $\chi' \in B_{T_t}(y_0)$,
\begin{equation}
	E[\rho(\mu_t, \mu')] = \sum_{j=1}^{\infty}2^{-j}E\left[\left|\int_{Y \times \hat{U}} q_j(y,u)\mu_t(dy,du) - \int_{Y \times \hat{U}} q_j(y,u)\chi'(dy,du)  \right|\right] \ge \delta
\end{equation}
implying, for large enough $N$,
\begin{equation}
	E[\rho(\mu_t, \mu')] = \sum_{j=1}^{N}2^{-j}E\left[\left|\int_{Y \times \hat{U}} q_j(y,u)\mu_t(dy,du) - \int_{Y \times \hat{U}} q_j(y,u)\chi'(dy,du)  \right|\right] \ge \frac{\delta}{2}.
\end{equation}
Hence, for any $\chi' \in B_{T_t}(y_0)$,
\begin{equation}
	\label{eq:thm_3_4_estimate}
	E\left[\sqrt{\sum_{j=1}^{N}\left|\int_{Y \times \hat{U}} q_j(y,u)\mu_t(dy,du) - \int_{Y \times \hat{U}} q_j(y,u)\chi'(dy,du)  \right|^2}\right] \ge \frac{c_N\delta}{2}, \quad c_N = \text{ constant.}
\end{equation}
Let $h(\cdot,\cdot)$ be the function introduced in Definition \ref{def:S-h-AC} with $j = N$ and define the sequence of elements $\kappa_t$ by the equation $\kappa_t \BYDEF \int_{Y \times \hat{U}}^{} h(y,u)\mu_t(dy,du)$. Note that, by (\ref{eq:V_h_G_representation}), $\kappa_t \in V_h$. Also, let $\kappa' \BYDEF \int_{Y \times \hat{U}}^{} h(y,u)\chi'(dy,du)$ for any $\chi' \in B_T(y_0)$ and note that $\kappa' \in V_h^r(T,y_0)$. The estimate (\ref{eq:thm_3_4_estimate}) is therefore equivalent to
\begin{equation}
	E[||\kappa_t - \kappa'||] \ge \frac{c_N \delta}{2}\quad \forall \ \kappa' \in V_h^r(T_t,y_0) \quad \iff \quad d^E(\kappa_t, V_h^r(T_t, y_0)) \ge \frac{c_N \delta}{2}
\end{equation}
which, on letting $t \rightarrow \infty$, contradicts (\ref{eq:V_h_and_E_V_h^r(T,y_0)_Edistance}), thus implying (\ref{eq:main_borkar_result_1}) holds. 
$\ \Box$

\section{Proof of Theorem \ref{thm:G_T_W_distance}} \label{sec:borkar_proofs_2}
\begin{Lemma}
	\label{lem:variance_Phi}
	
	Let $\pi \in \Pi$ be a control plan with the corresponding state-control trajectory $(y^{\pi, y_0}(\cdot), u^{\pi, y_0}(\cdot))$ given by (\ref{eq:associated_system}) - (\ref{eq:control_dependency}). Then, for any continuous function $\ph(\cdot): Y
	\rightarrow \mathbb{R}$, the following equalities are true:
	\begin{equation}
		\label{eq:exp_var_equation}
		E\left[\frac{1}{T}\sum_{t=0}^{T-1} \Phi(t)\right] = 0
	\end{equation}
	and
	\begin{equation}
		\label{eq:lim_exp_var}
		\lim_{T \rightarrow \infty}\textnormal{Var}\left(\frac{1}{T}\sum_{t=0}^{T-1} \Phi(t)\right) = 0,
	\end{equation}
	where 
	\begin{equation*}
		\Phi(t) = \ph(f(y^{\pi, y_0}(t),u^{\pi, y_0}(t),\sigma(t))) - \psi(y^{\pi,y_0}(t), u^{\pi,y_0}(t)),
	\end{equation*}
	with
	\begin{equation}\label{eq:psi-exp-phi}
		\psi(y, u) \BYDEF E[\ph (f(y,u,\sigma))].
	\end{equation}	
\end{Lemma}
{\bf Proof.}
Since $\sigma(t)$ is independent of the random elements $\sigma(0), \sigma(1), ..., \sigma(t-1)$ and since $y^{\pi,y_0}(t)$, $u^{\pi,y_0}(t)$ are functions of the latter,
$$
E[\ph (f(y^{\pi,y_0}(t), u^{\pi,y_0}(t),\sigma(t)))] = E[\psi(y^{\pi,y_0}(t), u^{\pi,y_0}(t))]\ \ \ \ \Rightarrow \ \ \ \ E[\Phi(t)] =0.
$$
The latter implies (\ref{eq:exp_var_equation}). From (\ref{eq:exp_var_equation}), it follows that
\begin{equation}
	\label{eq:variance_Phi}
	\textnormal{Var}\left(\frac{1}{T}\sum_{t=0}^{T-1} \Phi(t)\right)= E\left[\left(\frac{1}{T}\sum_{t=0}^{T-1} \Phi(t)\right)^2\right]=\frac{1}{T^2}\sum_{t=0}^{T-1}\sum_{t'=0}^{T-1} E[\Phi(t)\Phi(t')].
\end{equation}
Consider the terms of the sum on the right-hand side with $t \ne t'$ (assuming without loss of generality that $t > t'$). We have (see \cite[p. 216]{shiryaev1996probability})
\begin{equation*}
	E[\Phi(t)\Phi(t')] = E[E[\Phi(t)\Phi(t') \mid \Omega_t ]] = E[\Phi(t')E[\Phi(t) \mid \Omega_t ]],
\end{equation*}
where $\Omega_t$ is the $\sigma$-algebra generated by the random elements $\sigma(0), \sigma(1),...,\sigma(t-1)$. However (due to the fact that $\sigma(t)$ is independent of $\Omega_t$),
$$E[\ph (f(y^{\pi,y_0}(t), u^{\pi,y_0}(t),\sigma(t)))| \Omega_t]= \psi(y^{\pi,y_0}(t),u^{\pi, y_0}(t)) \quad \textit{a.s.}$$ 
Therefore,
\begin{equation*}
	E[\Phi(t) \mid \Omega_t ] =0  \quad \textit{a.s.}
\end{equation*}
Consequently, 
\begin{equation*}
	E[\Phi(t)\Phi(t')] = 0 \ \ \ \ \text{for $t \ne t'$,}
\end{equation*}
and, by (\ref{eq:variance_Phi}),
\begin{align*} 
	\textnormal{Var}\left(\frac{1}{T}\sum_{t=0}^{T-1} \Phi(t)\right)&=E\left[\frac{1}{T^2}\sum_{t=0}^{T-1}(\ph(f(y^{\pi, y_0}(t),u^{\pi, y_0}(t),\sigma(t))) - E [\ph(f(y^{\pi, y_0}(t),u^{\pi, y_0}(t),\sigma(t)))])^2\right] \\
	&\le  E\left[\frac{1}{T^2}\sum_{t=0}^{T-1}2\max_{y\in Y} |\ph(y)|^2\right] \\
	&\le E\left[\frac{2}{T}\max_{y\in Y} |\ph(y)|^2\right]
\end{align*}
which converges to zero as $T\rightarrow \infty$. Thus, (\ref{eq:lim_exp_var}) holds. $\ \Box$
\begin{Lemma}
	\label{lem:convergence_in_law}
	If a sequence of random variables $x_i$ is such that $E[x_i] \rightarrow E[x]$ and $\textnormal{Var}[x_i] \rightarrow 0$ as $i\rightarrow \infty$, then $x_i \rightarrow E[x]$ in law.
\end{Lemma}
{\bf Proof.}
By the Portmanteau Theorem (e.g. see \cite[p. ~254]{klenke2013probability}), $x_i \rightarrow E[x]$ in law if and only if for all bounded Lipschitz continuous functions $f(\cdot)$
\begin{equation}\label{eq:pmt_1}
	\lim_{i \rightarrow \infty}E[f(x_i)] = f(E[x]).
\end{equation}
Let $f(\cdot)$ be an arbitrary bounded Lipschitz continuous function with $L_f$ as the associated Lipschitz constant. Then, by the Cauchy-Schwartz Inequality,
\begin{align*}
	|E[f(x_i)] - f(E[x])| &\le E[|f(x_i) - f(E[x])|]\\
	&\le L_f E[|x_i - E[x]|]\\ 
	& \le L_f (E[|x_i - E[x_i]|] + |E[x_i] - E[x]|)\\ 
	& \le L_f \sqrt{E[|x_i - E[x_i]|^2]} + L_f|E[x_i] - E[x]|. 
\end{align*}
The first term on final line converges to zero since the variance of $x_i$ converges to zero by the premise. Similarly, the second term also converges to zero. Hence, (\ref{eq:pmt_1}) is satisfied, implying $x_i \rightarrow E[x]$ in law. $\ \Box$

\begin{Lemma}
	\label{lem:W_properties}
	The set $W$ (see (\ref{eq:s_stationary})) is presentable in the form
	\begin{equation}\label{eq:W_ctbl_dense_rep}
		W = \left\{\g \in \mathcal{P}(Y \times \hat{U}) : \int_{Y \times \hat{U}} \left(\ph(y) - E[\ph\big(f(y,u,\sigma)\big)]\right) \g(dy,du) = 0 \quad \ph \in \mathcal{D}(C(Y)) \right\}.
	\end{equation}
	where $\mathcal{D}(C(Y))$ is a countable dense set in $C(Y)$.
\end{Lemma}
{\bf Proof.} The proof is obvious.
$\ \Box$

\noindent
{\bf Proof of Theorem \ref{thm:G_T_W_distance}.}
Assume (\ref{eq:G_T_W_distance}) does not hold. Then, there exists a number $\delta > 0$, and initial condition $y_0$ such that, for some $T_i$, $\lim_{i \rightarrow \infty}T_i = \infty$, and some $\chi_i \in B_{T_i}(y_0)$,
\begin{equation}\label{eq:G_T_W_distance_contradiction}
	E\left[\inf_{\g' \in W}  \rho(\chi_i,\g')\right] \ge \delta \ \implies \ E\left[\inf_{\g' \in W}  \rho(\chi^*,\g')\right] \ge \delta,
\end{equation}
where it is assumed without loss of generality that $\chi_i \rightarrow \chi^*$ in law as $i \rightarrow \infty$.
From the fact that $\chi_i \in B_{T_i}(y_0)$, there exists a control plan $\pi_i \in \Pi$ with the corresponding state-control trajectory $(y^{\pi_i,y_0}(\cdot), u^{\pi_i,y_0}(\cdot))$ such that $\chi_i$ is the occupational measure of the pair $(y^{\pi_i,y_0}(\cdot), u^{\pi_i,y_0}(\cdot))$ on the interval $[0,T_i]$. Hence, for any $\ph \in \mathcal{D}(C(Y))$,
\begin{equation}\label{eq:G_T_W_distance_1}
	\begin{split}
		&\int_{Y \times \hat{U}} (\ph(y) - E [\ph(f(y,u,\sigma))]) \chi_i(dy,du) \\&= \frac{1}{T_i} \sum_{t=0}^{T_i - 1} \ph(y^{\pi_i,y_0}(t))  - \psi(y^{\pi_i,y_0}(t), u^{\pi_i,y_0}(t)) \\ 		
		&=\frac{1}{T_i} \sum_{t=0}^{T_i - 1} \ph(y^{\pi_i,y_0}(t)) - \ph(f(y^{\pi_i,y_0}(t),u^{\pi_i,y_0}(t),\sigma(t)))\\
		&+ \frac{1}{T_i}\sum_{t=0}^{T_i-1} (\ph(f(y^{\pi_i,y_0}(t),u^{\pi_i,y_0}(t),\sigma(t))) - \psi(y^{\pi_i,y_0}(t), u^{\pi_i,y_0}(t)),
	\end{split}
\end{equation}
where $\psi$ is as in (\ref{eq:psi-exp-phi}). The term in the second row of (\ref{eq:G_T_W_distance_1}) is a telescoping sum, hence,
\begin{equation}
	\label{eq:rho_mu_W_equation}
	\begin{split}
		&\int_{Y \times \hat{U}} (\ph(y) - E [\ph(f(y,u,\sigma))]) \chi_i(dy,du) \\
		&= \frac{1}{T_i}[\ph(y_0) - \ph(y^{\pi_i,y_0}(T_i))] +   \frac{1}{T_i}\sum_{t=0}^{T_i-1} \ph(f(y^{\pi_i,y_0}(t),u^{\pi_i,y_0}(t),\sigma(t))) - \psi(y^{\pi_i,y_0}(t), u^{\pi_i,y_0}(t)).
	\end{split}
\end{equation}
The first term in the second line of (\ref{eq:rho_mu_W_equation}) converges to zero as $i\rightarrow \infty$ by the compactness of $Y$. The expectation and the variance of the second term on the second line converges to zero by Lemma \ref{lem:variance_Phi}. 
Thus, by Lemma \ref{lem:convergence_in_law},
\begin{equation}
	\label{eq:law_converges_to_zero}\textbf{}
	\int_{Y \times \hat{U}} (\ph(y) - E[\ph(f(y,u,\sigma))]) \chi_i(dy,du) \rightarrow 0 \quad \textit{in law}
\end{equation}
for all $\ph \in \mathcal{D}(C(Y))$. Now, since $\chi_i \rightarrow \chi^*$ in law, by Skorohod's Theorem (e.g. see \cite[p. ~23]{borkar2012probability}, \cite[p. ~70]{billingsley2013convergence}) there exists random variables $\tilde{\chi}_i$ and $\tilde{\chi}^*$ on a common probability space where $\tilde{\chi}_i$ agrees in law with $\chi_i$; $\tilde{\chi}^*$ agrees in law with $\chi^*$, and 
\begin{equation}
	\lim\limits_{i\rightarrow \infty}\rho(\tilde{\chi}_i, \tilde{\chi}^*)=0 \quad \textit{ almost surely. }
\end{equation}
Hence, 
\begin{align*}
	&\int_{Y \times \hat{U}} ( \ph(y) - E [\ph(f(y,u,\sigma))] ) \tilde{\chi}^*(dy,du)  \\&= \lim\limits_{i\rightarrow \infty} \int_{Y \times \hat{U}} ( \ph(y) - E [\ph(f(y,u,\sigma))] ) \tilde{\chi}_i(dy,du) \quad \textit{a.s.}
\end{align*}
For each $i$, the term $\int_{Y \times \hat{U}} (\ph(y) - E [\ph(f(y,u,\sigma))]) \tilde{\chi}_i(dy,du)$ has the same distribution as the term on the left-hand side in (\ref{eq:law_converges_to_zero}), hence,
\begin{equation*}
	\int_{Y \times \hat{U}} ( \ph(y) - E [\ph(f(y,u,\sigma))] ) \tilde{\chi}^*(dy,du) = 0 \quad \textit{a.s.}
\end{equation*}
Since $\tilde{\chi}^*$ and $\chi^*$ agree in law, it follows that
\begin{equation}
	\label{eq:mu*_integral_zero}
	\int_{Y \times \hat{U}} ( \ph(y) - E [\ph(f(y,u,\sigma))] ) \chi^*(dy,du) = 0 \quad \textit{a.s.}
\end{equation}
Since $\mathcal{D}(C(Y)))$ is countable, the last expression is valid for all $\ph \in \mathcal{D}(C(Y)))$ outside a common zero probability set. That is, $\chi^* \in W$ a.s. (see (\ref{eq:W_ctbl_dense_rep})). Therefore,
\begin{equation}
	E\left[\inf_{\g' \in W}\rho(\chi^*,\g')\right]=0
\end{equation}	
which contradicts (\ref{eq:G_T_W_distance_contradiction}), completing the proof. $\ \Box$

\bibliographystyle{plain}
\bibliography{references}

\end{document}